\documentclass[final,5p,times,twocolumn]{elsarticle}

\makeatletter
\def\ps@pprintTitle{%
 \let\@oddhead\@empty
 \let\@evenhead\@empty
 \def\@oddfoot{\mbox{\emph{\small{}November 4, 2017}}}%
 \let\@evenfoot\@oddfoot}
\makeatother

\usepackage[utf8]{inputenc}
\usepackage[T1]{fontenc}
\usepackage[american]{babel}
\usepackage[babel=true]{microtype}
\usepackage{algorithm, algorithmic, amsmath, amssymb, siunitx, subcaption, tikz}
\usetikzlibrary{calc}
\usepackage{hyperref}
\hypersetup{%
   pdfauthor = {F. Dassi, L. Kamenski, P. Farrell, H. Si},
   pdftitle = {Tetrahedral mesh improvement using
      moving mesh smoothing, lazy searching flips,
      and RBF surface reconstruction},
   pdfsubject = {Preprint},
   pdfkeywords = {mesh improvement, mesh quality, edge flipping,
      mesh smoothing, moving mesh, radial basis functions}
}
\usepackage{xspace}
\usepackage{cleveref}
\usepackage{color}
\Crefname{ALC@unique}{Line}{Lines}
\crefname{algorithm}{Algorithm}{Algorithms}
\crefformat{section}{Section~#2#1#3}
\crefmultiformat{section}{Sections #2#1#3}{ and~#2#1#3}{, #2#1#3}{ and~#2#1#3}
\crefformat{figure}{Fig.~#2#1#3}
\Crefformat{figure}{Figure~#2#1#3}
\crefmultiformat{figs}{Figs #2#1#3}{ and~#2#1#3}{, #2#1#3}{ and~#2#1#3}
\crefrangeformat{figure}{Figs #3#1#4 to~#5#2#6}
\crefformat{equation}{Eq.~(#2#1#3)}
\Crefformat{equation}{Equation~(#2#1#3)}
\crefmultiformat{equation}{Eqs.~(#2#1#3)}{ and~(#2#1#3)}{, (#2#1#3)}{ and~(#2#1#3)}
\crefrangeformat{equation}{Eqs.~(#3#1#4) to~(#5#2#6)}

\newcommand{\Stellar}{\texttt{Stellar}}
\newcommand{\CGAL}{\texttt{CGAL}}

\newcommand{\MMG}{\texttt{mmg3d}}
\newcommand{\IR}{\mathbb{R}}                             
\DeclareMathOperator{\tr}{tr}
\providecommand{\abs}[1]{\lvert#1\rvert}
\providecommand{\Abs}[1]{\left\lvert#1\right\rvert}
\providecommand{\p}[2]{\frac{\partial{}#1}{\partial{}#2}}
\providecommand{\J}{\mathbb{J}}
\providecommand{\Th}{\mathcal{T}_h}
\providecommand{\Nv}{\#\mathcal{N}_h}
\providecommand{\Nt}{\#\mathcal{T}_h}
\providecommand{\bv}{\boldsymbol{v}}
\providecommand{\bx}{\boldsymbol{x}}
\providecommand{\bxi}{\hat{\boldsymbol{x}}}

\newcommand{\level}{\texttt{level}}

\newcommand{\imagePath}{.}

\newcommand{\plotResult}[9]{%
{\small{}
\begin{minipage}{0.23\textwidth}
\centering{}
\underline{#1}
\includegraphics[width=1.0\textwidth]{#9-#8-plot_dihedangles_normalized}
\vspace{-14.5em}
\begin{flushright}
$\theta_{\min,\mathcal{T}_h}=\ang{#2}\,\,\,$\\
$\theta_{\max,\mathcal{T}_h}=\ang{#3}\,\,\,$\\
$\mu_{\mathcal{T}_h}=\ang{#4}\,\,\,$\\
$\sigma_{\mathcal{T}_h}=#5\,\,\,\,\,$
\end{flushright}
\vspace{6.5em} 
\[
\#\mathcal{T}_h=\num{#7}
\]
\end{minipage}
}}

\newcommand{\plotResultDoubleSize}[9]{%
{\small{}
\begin{minipage}{0.46\textwidth}
\centering{}
\underline{#1}
\includegraphics[width=1.0\textwidth]{#9-#8-plot_dihedangles_normalized}
\vspace{-14.5em}
\begin{flushright}
$\theta_{\min,\mathcal{T}_h}=\ang{#2}\,\,\,$\\
$\theta_{\max,\mathcal{T}_h}=\ang{#3}\,\,\,$\\
$\mu_{\mathcal{T}_h}=\ang{#4}\,\,\,$\\
$\sigma_{\mathcal{T}_h}=#5\,\,\,\,\,$
\end{flushright}
\vspace{6.5em} 
\[
\#\mathcal{T}_h=\num{#7}
\]
\end{minipage}
}}

\newcommand{\plotAspectRatio}[3]{%
{\small{}
\begin{minipage}{0.4\textwidth}
\centering{}
\underline{#1}
\includegraphics[width=1.0\textwidth]{#3-#2}
\end{minipage}
}}

\newcommand{\plotTitle}{
\hfill{} \underline{Initial Mesh} \hfill{}\hfill{} \underline{New Method} \hfill{} \\}

\newcommand{\plotOnlyMesh}[3]{%
{\small{}
\begin{minipage}{0.23\textwidth}
\centering{}
\phantom{\underline{Initial Mesh}}
\includegraphics[width=0.95\textwidth]{#1-#2-#3}
\[
   \phantom{\#\mathcal{N}_h=\num{1},~~\#\mathcal{T}_h=\num{1}}
\]
\end{minipage}
}}

\newcommand{\plotOnlyMeshWithTitle}[4]{%
{\small{}
\begin{minipage}{0.23\textwidth}
\centering{}
\underline{#4}\\
\vspace{0.2cm}
\includegraphics[width=0.95\textwidth]{#1-#2-#3}
\[
   \phantom{\#\mathcal{N}_h=\num{1},~~\#\mathcal{T}_h=\num{1}}
\]
\end{minipage}
}}

\begin{document}


\begin{frontmatter}

\journal{}
\title{Tetrahedral mesh improvement using
   moving mesh smoothing,
   lazy searching flips,
   and~RBF surface reconstruction%
}

\author[bicocca]{Franco Dassi}
\ead{franco.dassi@unimib.it}

\author[wias]{Lennard Kamenski\corref{cor1}}
\ead{kamenski@wias-berlin.de}

\author[wias]{Patricio Farrell}
\ead{farrell@wias-berlin.de}

\author[wias]{Hang Si}
\ead{si@wias-berlin.de}

\cortext[cor1]{Corresponding author}

\address[bicocca]{Dipartimento di Matematica e Applicazioni,
   Università degli Studi di Milano-Bicocca,
   Via~Cozzi~53, 20125~Milano, Italy}

\address[wias]{Weierstrass Institute for Applied Analysis and Stochastics,
   Mohrenstr.~39, 10117~Berlin, Germany}


\begin{abstract}
Given a tetrahedral mesh and objective functionals measuring the mesh quality which take into account the shape, size, and orientation of the mesh elements, our aim is to improve the mesh quality as much as possible.
In this paper, we combine the \emph{moving mesh smoothing}, based on the integration of an ordinary differential equation coming from a given functional, with the \emph{lazy flip} technique, a reversible edge removal algorithm to modify the mesh connectivity.
Moreover, we utilize \emph{radial basis function} (RBF) surface reconstruction to improve tetrahedral meshes with curved boundary surfaces.
Numerical tests show that the combination of these techniques into a mesh improvement framework achieves results which are comparable and even better than the previously reported ones.
\end{abstract}


\end{frontmatter}


\section{Introduction}\label{sec:intro}

The key mesh improvement operations considered in this work are \emph{smoothing}, which moves the mesh vertices, \emph{flipping}, which changes the mesh topology without moving the mesh vertices, and a \emph{smooth boundary reconstruction}.
Previous work shows that the combination of smoothing and flipping achieves better results than if applied individually~\cite{FreitagOlliverGooch1997, Klingner:2007:ATM}.
In this paper, we combine the recently developed flipping and smoothing methods into one mesh improvement scheme and apply them in combination with a smooth boundary reconstruction via radial basis functions.

\emph{Mesh smoothing} improves the mesh quality by improving vertex locations, typically through Laplacian smoothing or some optimization-based algorithm.
Most commonly used mesh smoothing methods are Laplacian smoothing and its variants~\cite{Field1988,Lo1985}, where a vertex is moved to the geometric center of its neighboring vertices.
While economic, easy to implement, and often effective, Laplacian smoothing guarantees neither a mesh quality improvement nor mesh validity.
Alternatives are optimization-based methods that are effective with respect to certain mesh quality measures such as the ratio of the area to the sum of the squared edge lengths~\cite{Bank1994}, the ratio of the volume to a power of the sum of the squared face areas~\cite{SG1991},  the condition number of the Jacobian matrix of the affine mapping between the reference element and physical elements~\cite{FK02}, or various other measures~\cite{FreitagOlliverGooch1997,Knupp2000,Knupp1999,CTS98}.
Most of the optimization-based methods are local and sequential, combining Gauss-Seidel-type iterations with location optimization problems over each patch.
There is also a parallel algorithm that solves a sequence of independent subproblems~\cite{FJP1999}.

In our scheme, we employ the moving mesh PDE (MMPDE) method, defined as the gradient flow equation of a meshing functional (an objective functional in the context of optimization) to move the mesh continuously in time.
Such a functional is typically based on error estimation or physical and geometric considerations.
Here, we consider a functional based on the equidistribution and alignment conditions~\cite{Hua01b} and
employ the recently developed direct geometric discretization~\cite{Huang2015322} of the underlying meshing functional on simplicial meshes.
Compared to the aforementioned mesh smoothing methods, the considered method has several advantages:
it can be easily parallelized, it is based on a continuous functional for which the existence of minimizers is known, the functional controlling the mesh shape and size has a clear geometric meaning, and the nodal mesh velocities are given by a simple analytical matrix form.
Moreover, the smoothed mesh will stay valid if it was valid initially~\cite{HK2015}.

\emph{Flipping} is the most efficient way to locally improve the mesh quality and it has been extensively addressed in the literature~\cite{Joe1995-flips, FreitagOlliverGooch1997, GeorgeBoro2003-edgeflip, Klingner:2007:ATM}.
In the simplest case, the basic flip operations, such as 2-to-3, 3-to-2, and 4-to-4 flips, are applied as long as the mesh quality can be improved.
The more effective way is to combine several basic flip operations into one edge removal operation, which extends the 3-to-2 and 4-to-4 flips. This operation removes the common edge of $n \ge 3$ adjacent tetrahedra by replacing them with $m = 2n - 4$ new tetrahedra (the so-called $n$-to-$m$ flip).
There are at most $C_{n-2}$ possible variants to remove an edge by a $n$-to-$m$ flip, where $C_{n}=\frac{(2n)!}{(n+1)!\,n!}$ is the Catalan number.
If $n$ is small (e.g.,\ $n < 7$), one can enumerate all possible cases, compute the mesh quality for each case, and then pick the optimal one.
Another way is to use dynamic programming to find the optimal configuration.
However, the number of cases increases exponentially and finding the optimal solution with brute force is very time-consuming.

In this paper, we propose the so-called \emph{lazy searching flips}.
The key idea is to automatically explore sequences of flips to remove a given edge in the mesh.
If a flip sequence leads to a configuration which does not improve the mesh quality, the algorithm reverses this sequence and explores another one (see \cref{sub:edgeflip,fig:edgetoremove,fig:edgremoval,fig:treeflip}).
Once an improvement is found, the algorithms stops the search and returns without exploring the remaining possibilities.

When considering more arbitrary meshes (which may not be piecewise planar), we need to make sure that new nodes are added in a consistent way. To achieve this we use \emph{RBF surface reconstruction} as introduced in~\cite{Carr2001}. Radial basis functions are a very useful tool in the context of higher-dimensional interpolation as they dispense with the expensive generation of a mesh~\cite{Fornberg2015,Iske2004,Wendland2005}. Here, we will employ them to approximate the underlying continuous surface so that we can project nodes onto it as proposed in~\cite{Dassi2016,DasFarSi16}. This problem turns out to be very challenging for meshes with arbitrary boundary. Hence, we begin with a relatively simple mesh. For more complicated examples we first refine the boundary by using the RBF reconstruction and projection method and then keep the boundary nodes fixed while interior nodes may move.

In this paper, we provide a detailed numerical study of a combination of the MMPDE smoothing with the lazy searching flips and RBF surface reconstruction.
More specifically, we compare the results of the whole algorithm with \Stellar{}~\cite{Klingner:2007:ATM}, \CGAL{}~\cite{cgal} and \MMG{}~\cite{dobrzynski:hal-00681813}.
We also compare the lazy searching flips and the MMPDE smoothing with the flipping and smoothing procedures provided by \Stellar{}.

\section{The moving mesh PDE smoothing scheme}%
\label{sub:smoothing}

The key idea of this smoothing scheme is to move the mesh vertices via a moving mesh 
equation, which is formulated as the gradient system of an energy functional (the MMPDE 
approach).
Originally, the method was developed in the continuous setting~\cite{HRR94a,HR11}.
In this paper, we use its discrete form~\cite{Huang2015322, HK2015,HuaKamSi15},
for which the mesh vertex velocities are expressed in a simple, analytical matrix form, which makes the implementation more straightforward to parallelize.

\subsection{Moving mesh smoothing}

\begin{figure}[t]
   \centering{}
   \begin{tikzpicture}[scale = 1.0]
      %
      %
      \path ( 0.00,  0.00 ) coordinate (R1);
      \path ( 1.52,  0.00 ) coordinate (R2);
      \path ( 0.76,  1.32 ) coordinate (R3);
      \draw [fill = gray!12] (R1) -- (R2) -- (R3) -- cycle;
      \path ( 0.75,  0.15) coordinate (RK);
      \node [above] at (RK)  {$\hat{K}$};
      %
      %
      \path ( 3.50,  0.50 ) coordinate (M0);
      \path ( 4.00,  1.32 ) coordinate (M1);
      \path ( 6.00,  0.10 ) coordinate (M2);
      \draw [fill = gray!12] (M0) -- (M1) -- (M2) -- cycle;
      \path ( 3.80,  1.02 ) coordinate (MK);
      \node [below right] at (MK)  {$K$};
      %
      %
      \path ( 1.2, 1.02 ) coordinate (RMP);
      \path [->, line width = 1pt] (RMP)
         edge [bend left] node [above] {$F_K$} ($(M1) + (-0.40, -0.30)$);
      \path [->, line width = 1pt] ($(M0) + (0.0, -0.2)$)
         edge [bend left] node [above] {$F_K^{-1}$} ($(R2) + (0.15, 0.)$);
      %
         \node at (4.7, -0.2)  {$\J_K := {(F_K')}^{-1}$};
   \end{tikzpicture}
   \caption{Reference element $\hat{K}$, mesh element $K$,
      and the corresponding mappings $F_K$ and $F_K^{-1}$.\label{fig:notation}%
   }
\end{figure}
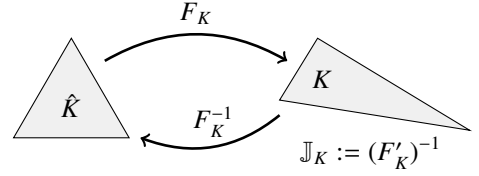

Consider a polygonal (polyhedral) domain $\Omega\subset \mathbb{R}^d$ with $d \ge 1$. Let $\Th$ denote the simplicial mesh as well as $\Nv$ and $\Nt$ the numbers of its vertices and elements, respectively.
Let $K$ be a generic mesh element and $\hat{K}$ the reference element taken as a regular simplex with volume $\abs{\hat{K}} = 1/\Nt$.
Further, let  $F_K'$ be the Jacobian matrix of the affine mapping $F_K \colon \hat{K} \to K$ from the reference element $\hat{K}$ to a mesh element $K$.
For notational simplicity, we denote the inverse of the Jacobian by $\J_K$, i.e.,\ $\J_K :=  {(F_K')}^{-1}$ (see \cref{fig:notation}).

Then, the mesh $\Th$ is uniform if and only if
\begin{equation}
   \abs{K} = \frac{\abs{\Omega}}{\Nt}
   \quad \text{and} \quad
   \frac{1}{d}\tr\left( \J_K^T \J_K \right) = \det{\left(\J_K^T \J_K\right)}^{\frac{1}{d}}
   \quad  \forall K \in \Th
   .
   \label{eq:equi:ali}
\end{equation}
The first condition requires all elements to have the same size and the second requires all elements 
to be shaped similarly to $\hat{K}$ (these conditions are the simplified versions of 
the equidistribution and alignment conditions~\cite{Hua06,HR11}).

The corresponding energy functional for which the minimization will result in a mesh satisfying \cref{eq:equi:ali} as closely as possible is
\begin{equation}
   I_h = \sum_{K} \Abs{K} \, G\left(\J_K, \det\J_K \right)
   \label{eq:Ih}
\end{equation}
with
\begin{equation}
   G(\mathbb{J}, \det\J) 
   = \theta  {\left(\tr\left(\J \J^T\right) \right)}^{\frac{dp}{2}}
         + \left(1 - 2 \theta\right) d^{\frac{dp}{2}} {(\det\J)}^p
   ,\label{eq:Ih:G}
\end{equation}
where $\theta \in (0,0.5]$ and $p > 1$ are dimensionless parameters (in \cref{sec:exe}, we use $\theta = 1/3$ and $p = 3/2$).
This is a specific choice and other meshing functionals are possible. The interested reader is referred to~\cite{HuaKamRus15} for a numerical comparison of meshing functionals for variational mesh adaptation.

In \cref{eq:Ih}, $I_h$ is a Riemann sum of a continuous functional for variational mesh adaptation based on equidistribution and alignment~\cite{Hua01b} and depends on the vertex coordinates $\bx_i$, $i= 1, \dotsc, \Nv$.
The corresponding vertex velocities $\bv_i$ for the mesh movement are defined as
\begin{equation}
   \bv_i := \frac{d \bx_i}{d t} = - {\left( \frac{\partial I_h}{\partial \bx_i}\right)}^T,
   \quad i = 1, \dotsc, \Nv,
   \label{eq:mmpde}
\end{equation}
where the derivatives $\frac{d \bx_i}{d t}$ are considered to be row vectors.

\subsection{Vertex velocities and the mesh movement}

The vertex velocities $\bv_i$ can be computed analytically~\cite[Eqs~(39) to~(41)]{Huang2015322} 
using scalar-by-matrix differentiation~\cite[Sect.~3.2]{Huang2015322}.
Denote the vertices of $K$ and $\hat{K}$ by $\bx_j^K$ and $\bxi_j$, $j = 0, \dotsc, d$, and define the element edge matrices as
\begin{align*}
   E_K &= [\bx_1^K - \bx_0^K, \dotsc, \bx_d^K - \bx_0^K]
   ,\\
   \hat{E} &= [\bxi_1 - \bxi_0, \dotsc, \bxi_d - \bxi_0]
   .
\end{align*}
Note, that $\hat{E} E_K^{-1} = \J_K$.
Then, the local mesh velocities are given element-wise~\cite[Eqs~(39) and~(41)]{Huang2015322} by
\begin{align}
   \begin{bmatrix}%
      {(\bv_1^K)}^T
      \\
      \vdots
      \\
      {(\bv_d^K)}^T
   \end{bmatrix}
    &= - G_K E_K^{-1}
      + E_K^{-1} \p{G_K}{\J} \hat{E} E_K^{-1}
      + \p{G_K}{\det\J} \frac{\det(\hat{E})}{\det(E_K)} E_K^{-1}
   ,\label{eq:vertex:velocity} \\
   {(\bv_0^K)}^T &= - \sum_{j=1}^d {(\bv_{j}^K)}^T
   ,\notag{}
\end{align}
where $G_K = G(\J_K, \det \J_K)$ and
\begin{align*}
   \p{G_K}{\J} &= \p{G}{\J}(\J_K) =  d p \theta  {\left( \tr(\J_K \J_K^T )\right )}^{\frac{d p}{2}-1}  \J_K^T
   ,\\
   \p{G_K}{\det\J} &= \p{G}{\det\J} (\det \J_K) = p (1-2\theta) d^{\frac{d p}{2}} {(\det\J_K)}^{p-1}
\end{align*}
are the derivatives of $G$ with respect to its first and second argument~\cite[Example~3.2]{Huang2015322} evaluated at $\J = \J_K$ and $\det(\J) = \det \J_K$.

The moving mesh equation~\eqref{eq:mmpde} becomes
\begin{equation}
   \frac{d \bx_i}{d t} 
      =  \sum\limits_{K \in \omega_i} \Abs{K} \bv_{i_K}^K,
   \quad i = 1, \dotsc, \Nv,
   \label{eq:mesh:eq}
\end{equation}
where $\omega_i$ is the patch of the vertex $\bx_i$ and $i_K$ is the local index of $\bx_i$ on $K$.

The moving mesh governed by \cref{eq:mesh:eq} will stay nonsingular if it is nonsingular initially: the minimum height and the minimum volume of the mesh elements will stay bounded from below by a positive number depending only on the initial mesh and the number of the elements~\cite{HK2015}.
This holds for the numerical integration of \cref{eq:mesh:eq} as well if the ODE solver has the property of monotonically decreasing energy~\cite{HK2015}.
For example, algebraically stable Runge-Kutta methods preserve this property under a mild step-size restriction~\cite{HaiLub14}.

During smoothing, we use the current vertex locations as the initial position and integrate \cref{eq:mesh:eq} for a time period (with the proper modification for the boundary vertices, see \cref{sec:velocity:adjustment}).
The connectivity is kept fixed during the smoothing step.
The time integration can be carried out for a given fixed time period or adaptively until the change of the energy functional \eqref{eq:Ih} is smaller than the prescribed absolute or relative tolerances, that is until
\begin{align*}
   | I_h(t_{n+1}) - I_h(t_{n}) | \le \varepsilon_{abs}
   \quad \text{or} \quad
   | I_h(t_{n+1}) - I_h(t_{n}) | \le \varepsilon_{rel} I_h(t_{n+1})
   .
\end{align*}

In our examples in \cref{sec:exe}, we use the explicit Runge-Kutta Dormand-Prince ODE solver~\cite{Dormand1980} and integrate until $t = 10$, which worked well with the provided examples.

\subsection{Velocity adjustment for the boundary vertices}
\label{sec:velocity:adjustment}

The velocities of the boundary vertices need to be modified.
If $\bx_i$ is a fixed boundary vertex, then its velocity is set to zero
Otherwise, $\bx_i$ is allowed to move along a boundary curve or a surface represented by the zero level set of a function $\phi$ and its velocity is modified so that its normal component along the curve (surface) is zero:
\[
   \nabla \phi (\bx_i) \cdot \p{\bx_i}{t} = 0.
\]

For the special case of a piecewise linear complex (PLC)~\cite{Miller1998} the velocity adjustment is straightforward:
\begin{align*}
   \text{facet vertices:}     &\quad \text{project the velocity onto the facet plane,}\\
   \text{segment vertices:}   &\quad \text{project the velocity onto the segment line,}\\
   \text{corner vertices:}    &\quad \text{set the velocity to zero.}
\end{align*}

For a general non-polygonal or non-polyhedral domain, a simple way to adjust the boundary vertices is to move the vertex and then project it onto the boundary to which it belongs, which proved to work well for simple surface geometries (see \cref{sec:curved:domains}).
However, for complicated geometries, this simple projection can fail and a more reliable approach is needed.

\section{Lazy searching flips}\label{sub:edgeflip}

In this section, we explain how to remove an edge and how to reverse the removal using flips.
In addition, we present the lazy searching algorithm which can be used to improve the quality of a mesh.

\subsection{Edge removal and its inverse}

A basic edge removal algorithm performs a sequence of elementary 2-to-3 and 3-to-2 flips~\cite{Si:2015:TDQ:2732672.2629697}. We extend this algorithm by allowing the flip sequence to be reversed. 
Our algorithm saves the flips online and it uses no additional memory.

Let $[a,b] \in {\cal T}_h$ be an edge with endpoints $a$ and $b$ and $A[0, \dotsc, n-1]$ be the array of $n \ge 3$ tetrahedra in ${\cal T}_h$ sharing $[a,b]$.
For simplicity, we assume that $[a,b]$ is an interior edge of ${\cal T}_h$, so that all 
tetrahedra in $A$ can be ordered cyclically such that the two tetrahedra 
$A[i]$ and $A[(i+1) \mod n]$ share a common face.
The index $i$ takes values in $\{0, 1, \dotsc, n-1\}$.
Throughout this section, additions involving indices will be modulo $n$.

Given such an array $A$ of $n$ tetrahedra, we want to find a sequence of flips that will 
remove the edge $[a,b]$. Moreover, we also want to be able to reverse this sequence
in order to return to the original state.

Our edge removal algorithm includes two subroutines
\begin{align*}
  [\texttt{done}, m] := &\texttt{flipnm}(A[0, \dotsc, n-1], \level)
  ,\\
  &\texttt{flipnm\_post}(A[0, \dotsc, n-1], m)
\end{align*}
with an array $A$ (of length $n$) of tetrahedra and an integer $\level$ defining the maximum recursive level as input.

\begin{description}
\item[\texttt{flipnm}]
executes ``forward'' flips to remove the edge $[a,b]$. It returns a Boolean value 
 indicating whether the edge is removed or not and an integer $m$ ($3 \le m \le n$).
If the edge is not removed ($\texttt{done} = \texttt{FALSE}$), $m$ indicates the current size of 
$A$ (initially, $m := n$).
\item[\texttt{flipnm\_post}]
must be called immediately after \texttt{flipnm}. It releases the memory allocated in \texttt{flipnm} and can perform ``backward'' flips to undo the flip sequence performed by \texttt{flipnm}.
\end{description}

\begin{figure*}[ht]
\begin{subfigure}[t]{0.18\textwidth}\centering{}
   \includegraphics[width=0.7\textwidth,clip]{\imagePath/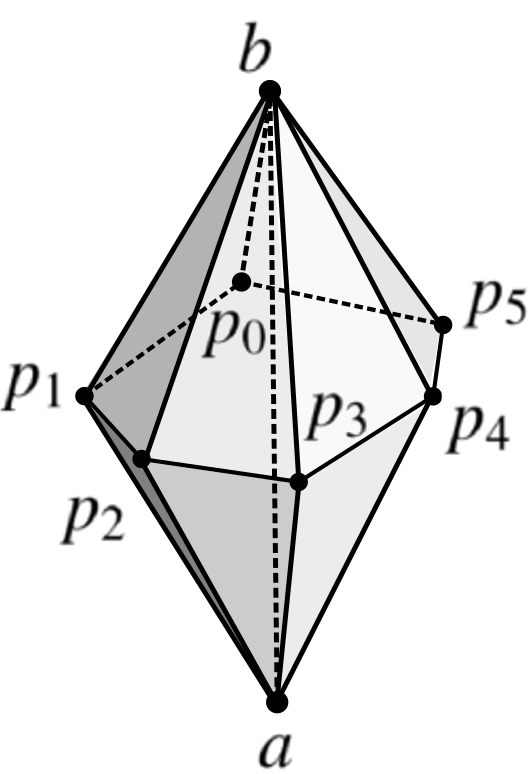}
   \caption{\label{fig:edgetoremove}%
      \footnotesize{}The initial state.%
   }
\end{subfigure}
\hfill{}
\begin{subfigure}[t]{0.80\textwidth}\centering{}
   \includegraphics[width=1.0\textwidth,clip]{\imagePath/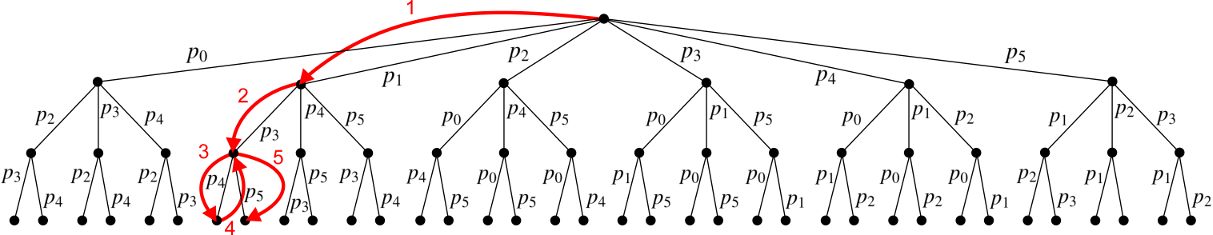}
   \caption{\label{fig:treeflip}%
      \footnotesize{}
      The lazy flip search tree for removing the edge $[a,b]$.
      $p_i$ identifies the face $[a,b,p_i]$ 
      which is flipped via a 2-to-3 flip. 
      The search path is highlighted with arrows.
   }
\end{subfigure}
\\[2ex]
\begin{subfigure}[t]{1.0\textwidth}\centering{}
   \begin{tabular}{cccc}
    \includegraphics[width=0.22\textwidth,clip]{\imagePath/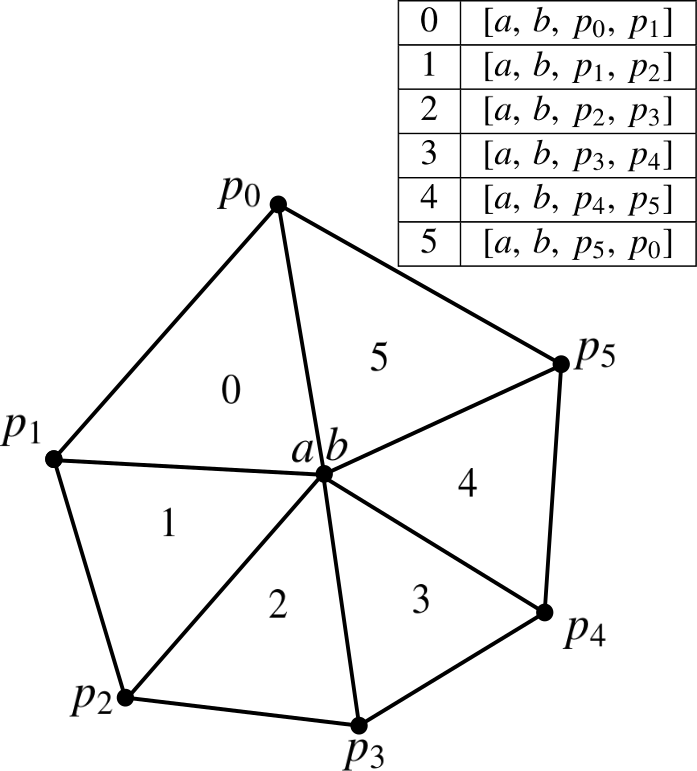}
   &\includegraphics[width=0.22\textwidth,clip]{\imagePath/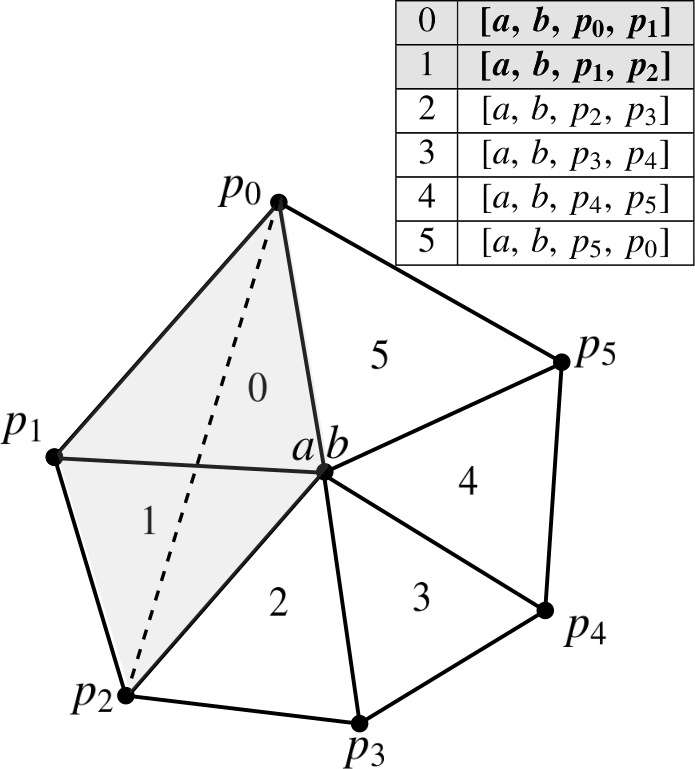}
   &\includegraphics[width=0.22\textwidth,clip]{\imagePath/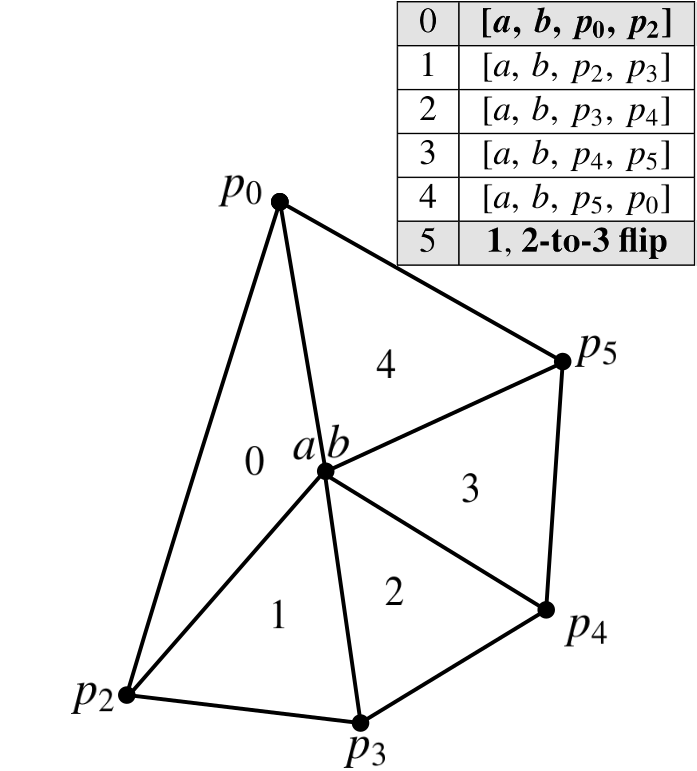}
   &\includegraphics[width=0.22\textwidth,clip]{\imagePath/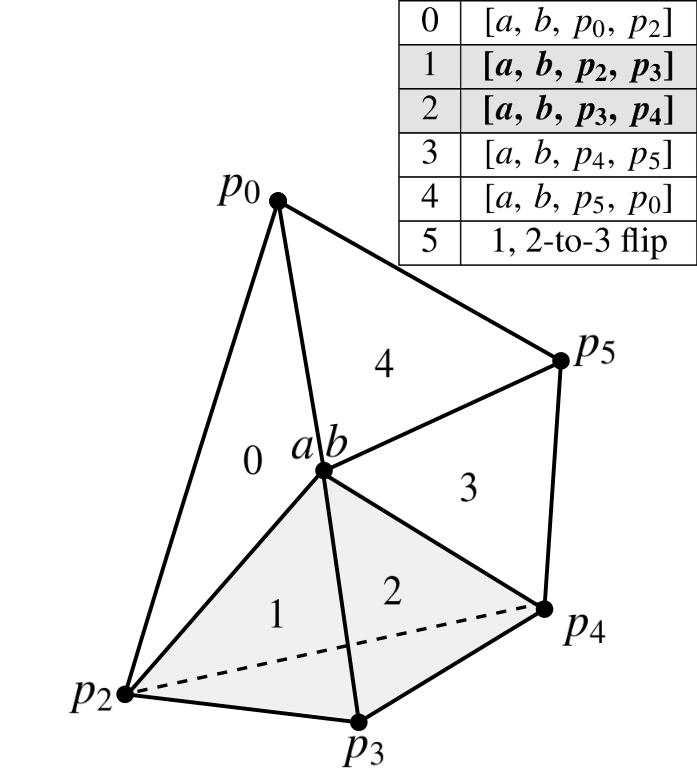}\\
   (1) & (2) & (3) & (4)
   \\[1ex]
    \includegraphics[width=0.22\textwidth,clip]{\imagePath/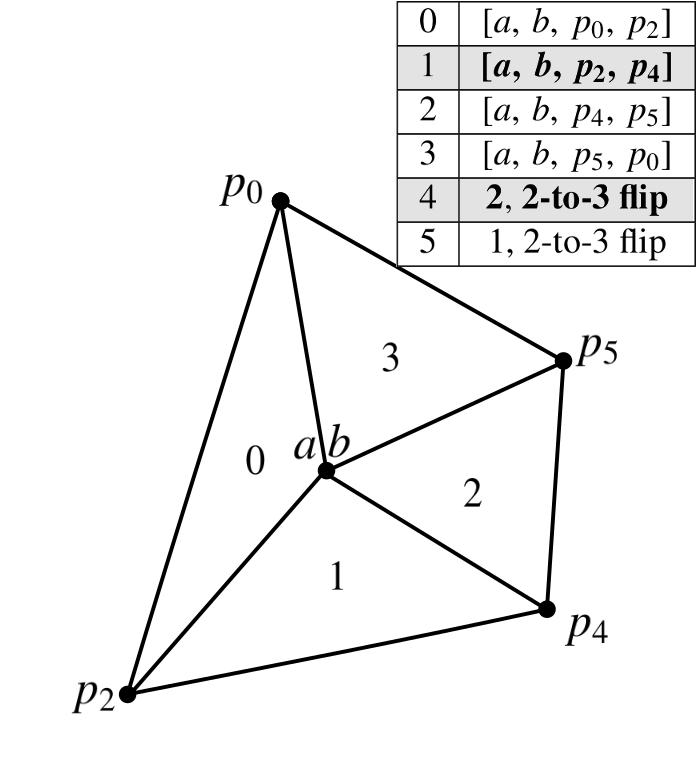}
   &\includegraphics[width=0.22\textwidth,clip]{\imagePath/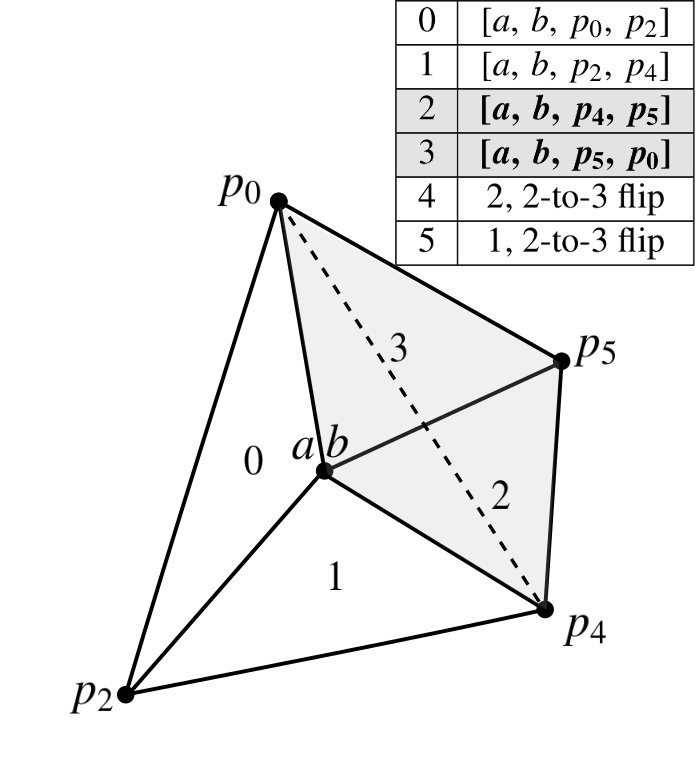}
   &\includegraphics[width=0.22\textwidth,clip]{\imagePath/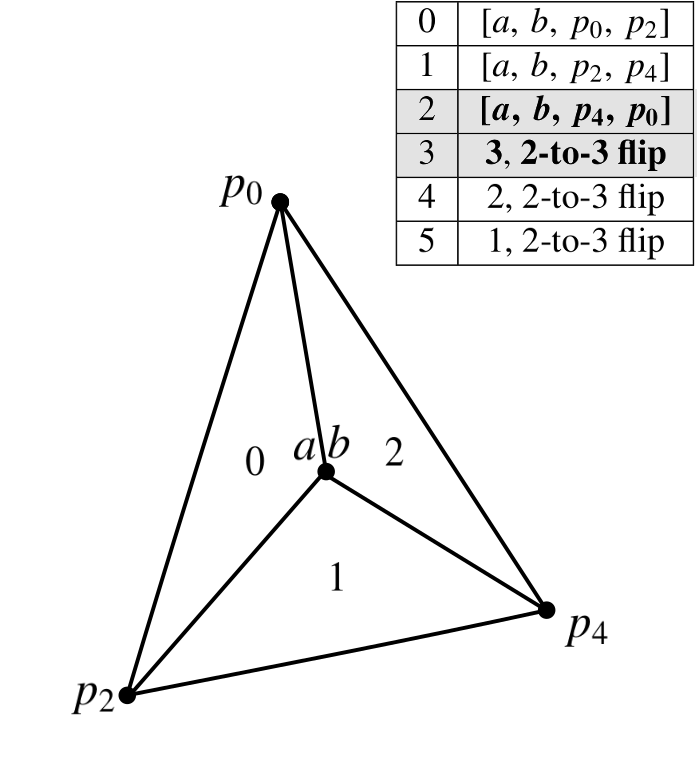}
   &\includegraphics[width=0.22\textwidth,clip]{\imagePath/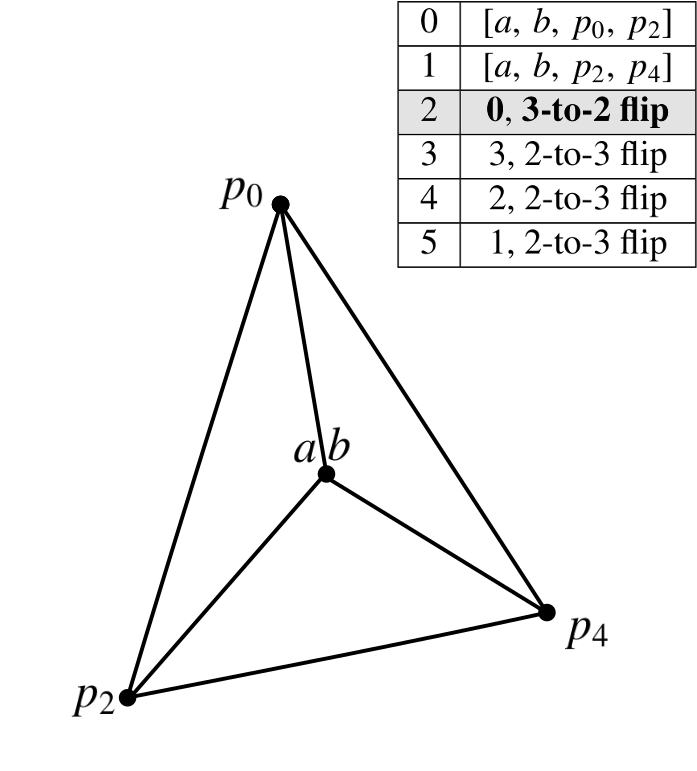}\\
   (5) & (6) & (7) & (8)
   \end{tabular}
   \caption{\label{fig:edgremoval}%
      \footnotesize{}
      The sequence of flips.
      The edge $[a,b]$ is represented by one vertex in the center (except (8)).
      A face $[a,b,p_i]$ is represented by an edge. 
      Arrays attached to each figure show the current content of $A$.
      (1) $n = 5$ tetrahedra share the edge $[a,b]$.
      In (2) and (3), $[a,b,p_1]$ is removed by a 2-to-3 flip.
      In (4) and (5), a 2-to-3 flip is performed on $[a,b,p_3]$.
      In (6) and (7), $[a,b,p_5]$ is removed by a 2-to-3 flip.
      In (8), the edge $[a,b]$ is removed by a 3-to-2 flip.
   }
\end{subfigure}
\caption{An example of an edge removal by a sequence of flips.}
\end{figure*}

The basic subroutine \texttt{flipnm}$(A[0, \dotsc, n-1],  \level)$ consists of the 
following three steps:

\begin{enumerate}[Step 1.]
\item 
Return $\texttt{done} = \texttt{TRUE}$ if $n = 3$ and \texttt{flip32} is possible for $[a,b]$
and $\texttt{done} = \texttt{FALSE}$ otherwise.

\item 
For each $i \in \{0, \dotsc, n-1\}$ try to remove the face $[a,b,p_i]$ by \texttt{flip23}.
If it is successfully flipped, reduce the size of $A$ by $1$.
Update $A[0, \dotsc, n-2]$ so that it contains the current set of tetrahedra sharing the edge $[a,b]$.
Reuse the last entry, $A[n-1]$, to store the information of this \texttt{flip23} (see Figure~\ref{fig:edgremoval}).
It then (recursively) calls $\texttt{flipnm}(A[0,\dotsc,n-2], \text{level})$.
When no face can be removed, go to Step~3.

\item If $\level > 0$, try to remove an edge adjacent to $[a,b]$ using \texttt{flipnm}.
For each $i \in \{0, \dotsc, n-1\}$, let $[x,y]$ be given by either edge $[a, p_i]$ or edge $[b, p_i]$. 
Initialize an array $B[0, \dotsc, n_1-1]$ of $n_1 \ge 3$ tetrahedra sharing $[x,y]$ and call \texttt{flipnm}$(B[0, \dotsc, n_1-1], \text{level}-1)$.
If $[x,y]$ is successfully removed, reduce $|A|$ by $1$.
Update $A[0, \dotsc, n-2]$ to contain the current set of tetrahedra sharing the edge $[a,b]$.
Reuse the last entry, $A[n-1]$, to store the information of this \texttt{flipnm}
and the address of the array $B$ (to be able to release the occupied memory later).
Then (recursively) call $\texttt{flipnm}(A[0, \dotsc, n-2])$.
Otherwise, if $[x,y]$ is not removed, call \texttt{flipnm\_post}$(B[0, \dotsc, n_1-1], m_1)$ to free the memory.
Return $\texttt{done} = \texttt{FALSE}$ if no edge can be removed.
\end{enumerate}

Since \texttt{flipnm} is called recursively, not every face and edge should be flipped in Steps 2 and 3.
In particular, if $B$ is allocated, i.e., \texttt{flipnm} is called recursively, we skip flipping faces and edges belonging to the tetrahedra in $A \cap B$.

In the simplest case, that is, ignoring the option to reverse the flips, 
\texttt{flipnm\_post}$(A[0, \dotsc, n-1], m)$ simply walks through the array $A$ from $A[m]$ to $A[n-1]$ and checks if a $\texttt{flipnm}$ flip has been saved. If so, the saved array address $B$ is extracted and its memory is released.

In Step~2 there are at most ${n \choose n - 3} / (n - 3)!$ different flip 
sequences, depending on the specific choice of faces in $A$.
Each individual flip sequence is equivalent to a sequence of the $n$ vertices (apexes) in the link of the edge $[a,b]$. 
We reuse the entries of $A$ to store each flip sequence. 
After a 2-to-3 flip, the number of the tetrahedra in array $A$ is reduced by one (two tetrahedra out, one terahedron in), since only one of the three new tetrahedra contains the edge $[a,b]$.  The remaining tetrahedra are shifted by one in the list, so that the last entry, $A[n-1]$, can be used to store this flip (cf.~\cref{fig:edgremoval}).
In particular, the following information is saved:
\begin{itemize}
\item a flag indicating a 2-to-3 flip;
\item the original position $i$, meaning that the face $[a,b,p_i]$ is flipped.
\end{itemize}
Both is compressed and stored in the entry $A[n-1]$ (note that a flag requires just a few bits of space).
This stored data allows us to perform the reversal of a 2-to-3 flip as follows:
\begin{itemize}
\item use $A$ and the position $i$ with
\[
      A[i-1] = [a,b,p_{i-1},p_{i+1}]
\]
to locate the three tetrahedra sharing the edge $[p_{i-1}, p_{i+1}]$: $[p_{i-1}, p_{i+1}, a,b]$, $[p_{i-1}, p_{i+1}, b, p_{i}]$, and  $[p_{i-1}, p_{i+1}, p_{i}, a]$;
\item perform a 3-to-2 flip on these three tetrahedra;
\item insert two new tetrahedra into the array $A$: 
\begin{align*}
   A[i-1] &= [a,b,p_i,p_{i-1}]
   ,\\
   A[i] &= [a,b,p_i,p_{i+1}]
.
\end{align*}
\end{itemize}

In Step 3, if the selected edge $[a,p_i]$ is removed, the sequence of flips to remove $[a,p_i]$ is stored in $B$. We then use the last entry $A[n-1]$ to store this sequence of flips.
In particular, the following information is saved:
\begin{itemize}
   \item a flag indicating that this entry stores the flip sequence to remove the edge $[a, p_i]$;
\item the original position $i$, i.e., the edge $[a,p_i]$ is flipped;
\item the address of the array $B$ which stores the flip sequence.
\end{itemize}
This information allows us to reverse this sequence of flips exactly.

\subsection{Lazy searching flips}

During the mesh improvement process we perform flips to improve the mesh quality.
Let us consider the case when it becomes necessary to remove an edge.
The maximum possible number of flips for an edge removal is the Catalan number $C_{n-2}$ ($n$ is the size of $A$).  Hence, the direct search for the optimal solution is only meaningful if $n$ is very small. 
In most situations, an edge may not be flipped if we restrict ourselves to adjacent faces 
of the edge.  
Our strategy is to search and perform the flips as long as they improve the current mesh quality. Our lazy searching scheme is not restricted by the number $n$ and can be extended to adjacent edges. 

The lazy searching flip scheme is like a walk in a $k$-ary search tree (a rooted 
tree with at most $k$ children at each node, see \cref{fig:treeflip}).
The root represents the edge $[a,b]$ to be flipped and each of the tree nodes represents 
either an adjacent face $[a,b,p_i]$ or an adjacent edge $[a,p_i]$ or $[b,p_i]$ of $[a,b]$.  
The edges of the tree represent our search paths. In particular, the directed edge from 
$\level$ $l$ to $l+1$ represents either a \texttt{flip23} or a \texttt{flipnm}, and the 
reversed edge represents the inverse operation. The tree depth is given by the parameter $\level$.

At $\level > 0$, in order to to decide if an adjacent face $[a,b,p_i]$ should be flipped, we 
check if $[a,b,p_i]$ is flippable and make sure that this flip improves the local 
mesh quality.  Note that we need to check only two of the three new tetrahedra: 
$[a,p_{i-1},p_i,p_{i+1}]$ and $[b,p_{i-1},p_i,p_{i+1}]$. The tetrahedron $[a,b,p_{i-1}, 
p_{i+1}]$ will be involved in the later flips, and will be flipped if the edge $[a,b]$ is 
flipped.

Once an improvement is found, the algorithm moves on to the next edge without exploring other possibilities.

\section{Radial basis functions to handle curved boundaries}

We describe in this section how to project the mesh on a smooth surface 
in order to deal with curved boundaries.
We achieve this with the help of radial basis functions (RBFs), see~\cite{Fornberg2015,Iske2004,Wendland2005}.

\subsection{Basic concepts and examples}

Let $\mathcal{P}_m(\IR^d)$ denote the space of $d$ variate polynomials 
with absolute degree at most $m$ 
and dimension $q:=\dim \mathcal{P}_m(\IR^d)=\binom{m-1+d}{d}$. For 
a basis $p_1,\hdots,p_q$ of this space, 
define the $M\times q$ polynomial matrix $P_X$ through its $ij^{th}$ entry,
\[
	p_{ij} = p_i(\mathbf{x}_j)\,,
\]
where $\mathbf{x}_j\in X$ and $X=\{\mathbf{x}_1,\hdots,\mathbf{x}_M\} \subseteq \IR^d$ 
denotes a data set. The function $\Phi$ is called 
\textit{conditionally positive definite of order m}
if the quadratic form
\begin{align*}
   {\mathbf{c}}^T A_{\Phi,X} \mathbf{c}
\end{align*}
for the distance matrix $A_{\Phi,X}$ with its $ij^{th}$ entry defined by
\[
   {(A_{\Phi,X})}_{ij} = \Phi(\mathbf{x}_i-\mathbf{x}_j)\,,
\]
is positive for all data sets $X$ and 
for all $\mathbf{c}\in\IR^M\setminus\{\mathbf{0}\}$ which 
additionally satisfy the constraint $P^T_X \mathbf{c} = \mathbf{0}$.

Conditionally positive functions of order $m$ are also conditionally positive definite
for any order higher than $m$.
Hence, the order shall denote the smallest positive integer $m$.
A conditionally positive 
definite function of order $m=0$ is called \textit{positive definite}.

\begin{figure}[ht]\centering{}
   \includegraphics[width=0.6\linewidth]{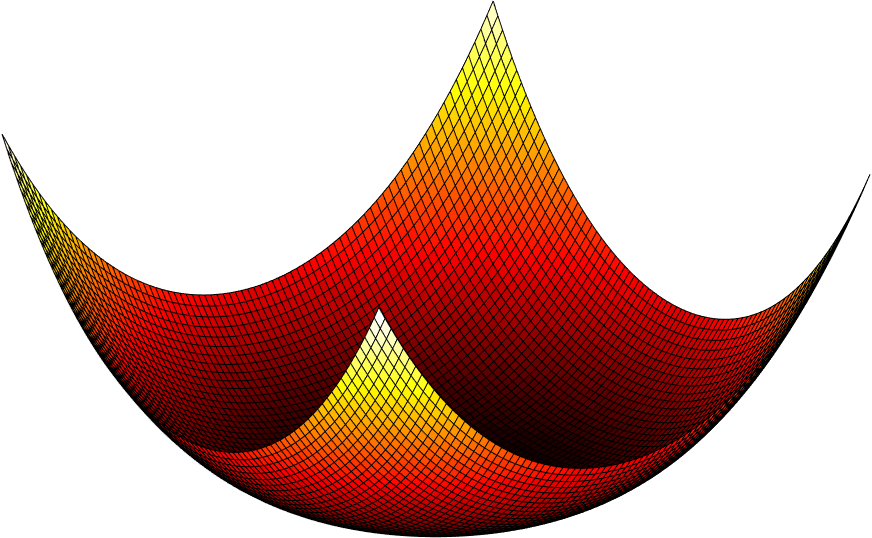}
   \caption{The polyharmonic spline $\|\mathbf{x}\|^3$.\label{fig:tps}}
\end{figure}

One speaks of radial basis functions if one 
additionally assumes
that $\Phi$ is a radial function, i.e.,\ there exists a function
$\phi\colon\IR_{\geq 0}\to \IR$
such that $\Phi(\mathbf{x})=\phi(\|\mathbf{x}\|)$. 
Common examples of RBFs include:
\begin{align*}
 \text{Gaussian:}               &\quad e^{-\|\mathbf{x}\|^2},\\
 \text{Multiquadric:}           &\quad \sqrt{1 + \|\mathbf{x}\|^2},\\
 \text{Inverse Multiquadric:}   &\quad 1 / \sqrt{1+\|\mathbf{x}\|^2},\\
 \text{Polyharmarmonic Spline:} &\quad \|\mathbf{x}\|^3.
\end{align*}
For the numerical examples in this paper, 
we exclusively use the polyharmonic spline $\|\mathbf{x}\|^3$ (\cref{fig:tps})
which is conditionally positive of order 2.

We assume now that the interpolant $s\colon\IR^d\to \IR$ is given by a linear combination
of translated radial basis function, augmented by a polynomial part, i.\,e.
\begin{equation}
	s(\mathbf{x}) = \sum_{j=1}^M \alpha_j\Phi(\mathbf{x}-\mathbf{x}_j) + 
\sum_{k=1}^q \beta_k p_k(\mathbf{x}).
\label{eq:cpdinterpolant}
\end{equation}
Thus, we have $M+q$ unknown coefficients, $M$ of which 
are determined from the interpolation conditions and 
$q$ conditions from requiring that $P^T_X \mathbf{c} = \mathbf{0}$.
For positive definite
functions, the linear system is positive definite by construction. 
Hence the coefficients can be determined uniquely.
It is also not difficult 
to verify that the interpolation and polynomial constraint conditions for conditionally 
positive definite functions lead to a uniquely solvable system,
see~\cite[Theorem 8.21]{Wendland2005} for details. In the case of conditionally 
positive definite functions, it is known that at least $M-q$ eigenvalues of 
the matrix $A_{\Phi,X}$ are positive~\cite[Section 8.1]{Wendland2005}.

\subsection{Surface reconstruction with RBFs}\label{sub:rfbSurf}

We will assume that the surface $\Gamma$ is given implicitly by the zero level set of some function
$F\colon\Omega \subseteq \mathbb{R}^3 \to \mathbb{R}$, i.\,e.
\begin{equation}
\Gamma = \left\{{(x,y,z)}^T \in \Omega \mid F(x,y,z)=0\right\}\, ,
\label{eq:implicitSurface}
\end{equation}
for some bounded domain $\Omega$.

We cannot simply assume that the target function (which we wish to interpolate) 
is the zero level set of the function $F$ since the right-hand 
side of the linear system one needs to solve would vanish which in turn implies 
that the coefficients vanish as well. Carr et al.~\cite{Carr2001} therefore made the additional 
assumption that the normal vectors are known. Then one can also prescribe 
on-surface and off-surface points. 
Assume that the points on 
the surface are denoted with $X=\{\mathbf{x}_1,\hdots,\mathbf{x}_N\}$ and the 
corresponding normal vectors with $M=\{\mathbf{n}_1,\hdots,\mathbf{n}_N\}$. We 
define the surface interpolation problem 
\begin{equation}
\begin{aligned}
	s(\mathbf{x}_i) &= F(\mathbf{x}_i) = 0,  
	&  1&\leq i\leq N \\
	s(\mathbf{x}_i+\varepsilon\mathbf{n}_i) &= 
F(\mathbf{x}_i+\varepsilon\mathbf{n}_i) = \varepsilon,  &   N+1&\leq i\leq 2N 
\end{aligned}
\label{eq:interpolCondition2}
\end{equation}
for some parameter $\varepsilon>0$. Since the right-hand side of the linear system does 
not vanish anymore, we find a nontrivial solution.
Recently, this surface interpolation technique was combined with the 
higher dimensional embedding technique~\cite{Dassi2016,DasFarSi16} to construct curvature-aligned anisotropic surface meshes. 
In this context the data set $X$ corresponds to the vertices of the mesh. 

\subsection{Projection onto the reconstructed surface}

\begin{figure}[!htb]\centering{}
\hfill{}
\begin{subfigure}[t]{0.29\linewidth}
   \includegraphics[width=1.0\linewidth]{\imagePath/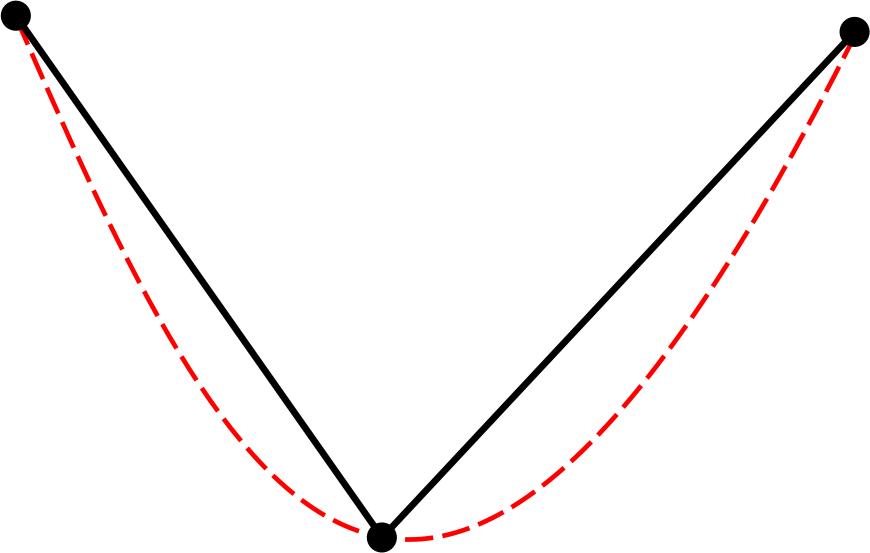}
   \subcaption{}
\end{subfigure}
\hfill{}
\begin{subfigure}[t]{0.29\linewidth}\centering{}
   \includegraphics[width=1.0\linewidth]{\imagePath/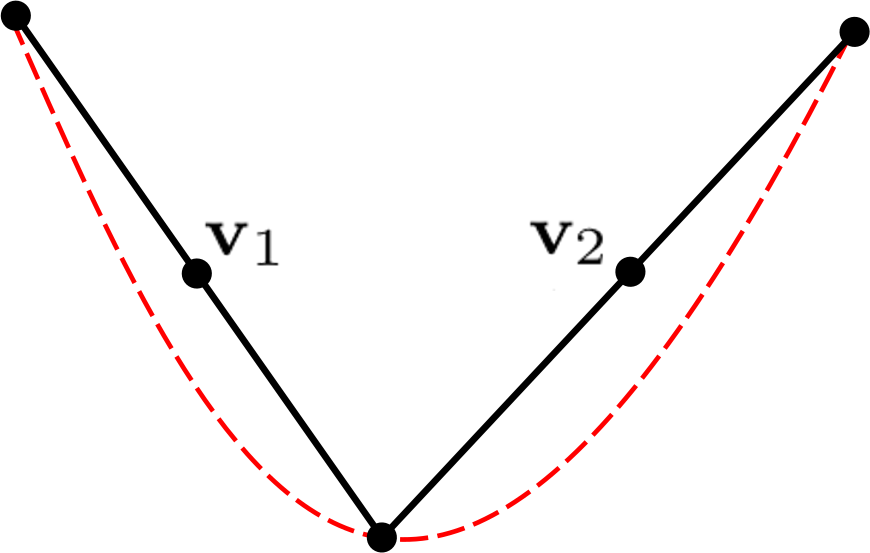}
   \subcaption{}
\end{subfigure}
\hfill{}
\begin{subfigure}[t]{0.29\linewidth}\centering{}
   \includegraphics[width=1.0\linewidth]{\imagePath/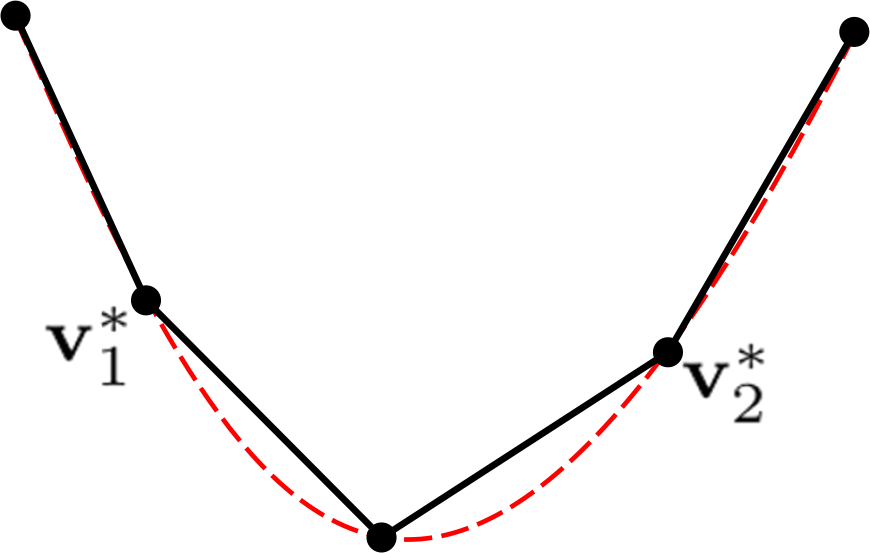}
   \subcaption{}
\end{subfigure}
\hfill{}
\caption{%
Edge splitting.
In (a), a coarse input mesh (solid line) approximates the reconstructed curve (dashed line).
In (b), the edges are halved, the midpoints $\mathbf{v}_1$ and $\mathbf{v}_2$ are not on the curve.
In (c), new points $\mathbf{v}_1^*$ and $\mathbf{v}_2^*$ are projected onto the curve using a steepest descent method.}\label{fig:projCount}
\end{figure}

There are two important parts of the projection algorithm:
edge splitting and edge contraction.
If we split an edge or move a point during
smoothing, we project the resulting point onto the RBF surface reconstruction (\cref{fig:projCount}).
\\
The projection itself is realized with ideas from~\cite{Hartmann1999355}.
This procedure is a combination of orthogonal projections on tangent planes
as well as tangent parabolas.
It requires only first order derivatives and uses a steepest descent method.
The combination of this projection method with RBF surface reconstruction has also been
discussed in~\cite{Dassi2016,DasFarSi16}.

\begin{figure}[!htb]\centering{}
\hfill{}
\begin{subfigure}[t]{0.29\linewidth}
   \includegraphics[width=1.0\linewidth]{\imagePath/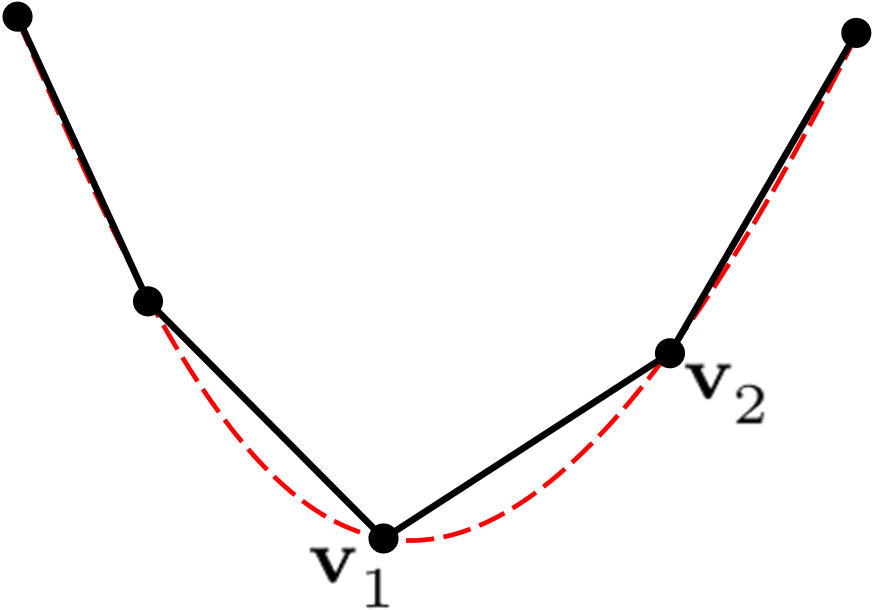}
   \subcaption{}
\end{subfigure}
\hfill{}
\begin{subfigure}[t]{0.29\linewidth}\centering{}
   \includegraphics[width=1.0\linewidth]{\imagePath/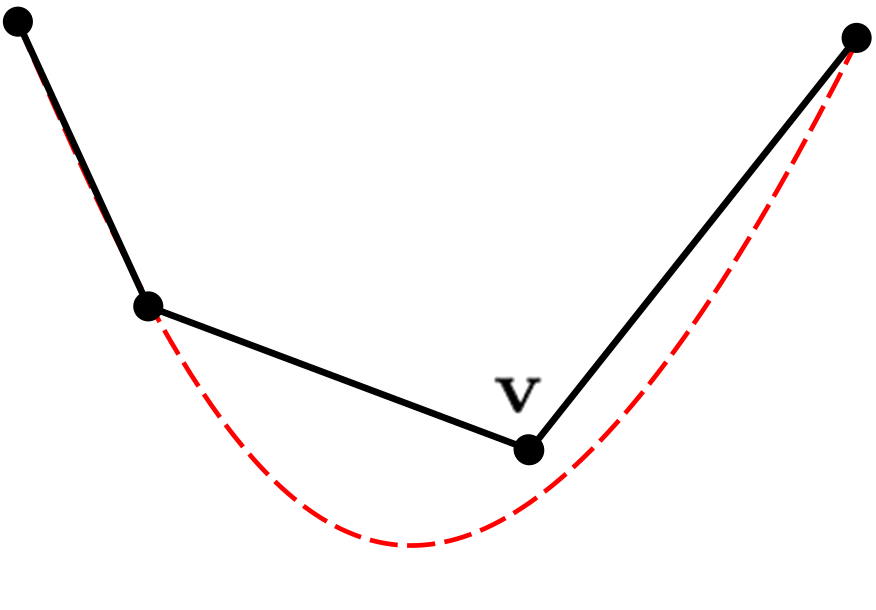}
   \subcaption{}
\end{subfigure}
\hfill{}
\begin{subfigure}[t]{0.29\linewidth}\centering{}
   \includegraphics[width=1.0\linewidth]{\imagePath/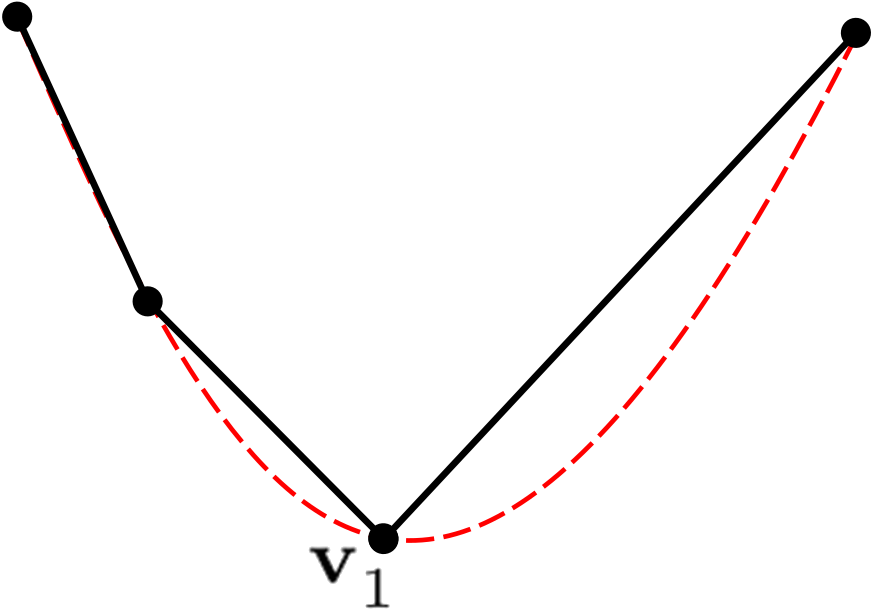}
   \subcaption{}
\end{subfigure}
\hfill{}
\caption{%
Edge contraction.
In (a), a fine input mesh (solid line) approximates the reconstructed curve (dashed line).
In (b), the edge $\overline{\mathbf{v}_1\mathbf{v}_2}$ is contracted into its midpoint $\mathbf{v}$, which is not on the curve.
In (c), the edge $\overline{\mathbf{v}_1\mathbf{v}_2}$ in contracted into its end point $\mathbf{v}_1$, which, by construction, is on the curve.}\label{fig:collIssue}
\end{figure}

When it becomes necessary to  contract an edge, we contract it into one of its endpoints (\cref{fig:collIssue}).

\section{Mesh improvement strategy}\label{sec:improve}

The goal of the proposed algorithm is to obtain a new isotropic mesh whose elements come ``as close as possible'' to the equilateral one.  To achieve this goal, we combine the local and global mesh operations described in \cref{sub:smoothing,sub:edgeflip}.

\subsection{Mesh quality}\label{sub:qual}

To say ``as close as possible to an equilateral tetrahedron'' is somewhat vague from a mathematical point of view.
To have a more precise criterion, the majority of the mesh improvement algorithms define a computable quantity $q(K)$ which quantifies how far a tetrahedron $K$ is from being equilateral~\cite{Klingner:2007:ATM,dobrzynski:hal-00681813,gmsh,MR3043640,Shewchuk02whatis,NME:NME804}.
Here, we take into account the following two:

\begin{description}
\item[Aspect Ratio:]
This is one of the most classical ways to evaluate the quality of a tetrahedron.
It is defined as
\begin{equation}\label{eqn:aspectRatio}
   q_{ar}(K) := \sqrt{\frac{2}{3}} \frac{L}{h} ,
\end{equation}
where $L$ is the longest edge and $h$ is the shortest altitude of the tetrahedron $K$.
By construction, $q_{ar}(K)\geq 1$ and an equilateral tetrahedron is characterized by $q_{ar}(K)=1$.
\item[Min-max Dihedral Angle:]
For each tetrahedron $K$ we consider both the minimal and the maximal dihedral angles $\theta_{\min,K}$ and $\theta_{\max,K}$.
An equilateral tetrahedron has $\theta_{\min,K}=\theta_{\max,K}=\arccos{(1/3)}\approx \ang{70.56}$.
Applying an operation that increases $\theta_{\min,K}$ or decreases $\theta_{\max,K}$ of a given tetrahedron $K$ makes $K$ ``closer'' to the equilateral shape.
Note that this is not a classical quality measure since we associate two quantities with each tetrahedron, which is one of the novel aspects of the proposed mesh improvement procedure.
\end{description}

These two quality measures refer to a single tetrahedron $K$ of the mesh. 
However, the design of our mesh improvement scheme requires a quality measure
for the \emph{whole} mesh as a stopping criterion.
To estimate the quality of the whole mesh, we define the global parameter
\begin{equation}\label{eqn:qualmesh}
   Q(\Th):= \min_{K\in\Th}\left( \theta_{\min,K} \right)
   .
\end{equation}
If we consider a target dihedral angle $\theta_{\lim}$ and 
obtain a mesh $\Th$ with $Q(\Th)>\theta_{\lim}$, 
then \emph{all} dihedral angles are guaranteed to be greater 
than $\theta_{\lim}$.

\subsection{The scheme}\label{sub:scheme}

The inputs for the mesh improvement algorithm are a tetrahedral mesh $\Th^{ini}$ of a PLC and a target minimum angle $\theta_{\lim}$.
The output is a mesh $\Th^{fin}$ where each element has a minimum dihedral angle greater than $\theta_{\lim}$.

\begin{algorithm*}[p]
  \textsc{Improve}$(\Th^{ini}$,\,$\theta_{\lim})$
  \begin{algorithmic}[1]
  \REPEAT{}\label{pt:bigIni}
    \fbox{\parbox{0.82\linewidth}{
    \REPEAT{}\label{pt:allOpIni} 
      \REPEAT{}
	\REPEAT{}
	  \fbox{\parbox{0.75\linewidth}{
	  \REPEAT{}\label{pt:smooth&flipStart}
     \STATE{} \emph{MMPDE-based smoothing}
     \STATE{} \emph{RBF surface reconstruction}
     \STATE{} \emph{lazy flips}
	  \UNTIL{no point is moved or no flip is done or 
	  $Q(\Th)\geq\theta_{\lim}$}\label{pt:smooth&flipEnd}
	  }}$\to$ smooth and flip
	  \STATE{} remove the edges $l_{\textbf{e}} < 0.5\,l_{\textrm{ave}}$\label{pt:coll}
	\STATE{} \emph{lazy flips}\label{pt:allcoll}
	\UNTIL{no edge is contracted or $Q(\Th)\geq\theta_{\lim}$}
	\STATE{} split the edges $l_{\textbf{e}} > 1.5\,l_{\textrm{ave}}$\label{pt:split}
   \STATE{} \emph{RBF surface reconstruction}
	\STATE{} \emph{lazy flips}\label{pt:allsplit}
      \UNTIL{no edge is split or $Q(\Th)>\theta_{\lim}$}
      \STATE{} split the tetrahedra $K$ such that $\theta_{\min,K}<\theta_{\lim}$\label{pt:random}
      \STATE{} \emph{lazy flips}\label{pt:allrandom}
    \UNTIL{no tetrahedron is removed or 
    $Q(\Th)>\theta_{\lim}$}\label{pt:allOpEnd}
    }}$\to$ main loop
    \STATE{} change the flip criterion for the \emph{lazy flips}\label{pt:changeLazyFlip}
   \UNTIL{no operation is done in the main loop or 
   $Q(\Th)>\theta_{\lim}$\label{pt:bigEnd}}
  \end{algorithmic}
  \caption{The proposed mesh improvement scheme.\label{alg:allProc}}
\end{algorithm*}

The scheme is presented in \cref{alg:allProc} and consists of five nested ``\textbf{repeat \ldots{} until}'' loops, whose stopping criterion depends on the operations done inside the loop and $Q(\Th)$.
We apply the MMPDE smoothing and the \emph{lazy flip} in the most internal loop (\crefrange{pt:smooth&flipStart}{pt:smooth&flipEnd}).
The \emph{lazy flip} is also exploited in the outer loops both on the whole mesh (\cref{pt:allcoll,pt:allsplit,pt:allrandom}) and on the tetrahedra involved in the local operations (\cref{pt:coll,pt:split,pt:random}).

It is possible to consider several flipping criteria for the \emph{lazy flip}, which makes the design of the scheme flexible.
We exploit this feature by using two objective functionals and changing the flipping criterion in each iteration of the outer loop (\cref{pt:changeLazyFlip}) by
\begin{enumerate}
   \item maximizing $\theta_{\min,K}$ and minimizing $\theta_{\max,K}$ (simultaneously),
   \item minimizing the aspect ratio.
\end{enumerate}
The stopping criterion is always based on the minimal dihedral angle, $Q(\Th)$, and the number of operations done.

After a number of iterations both the flipping and the smoothing procedure can stagnate, i.e.,\  the mesh $\Th$ converges to a fixed configuration where neither flips nor smoothing can improve the quality of the mesh.
Unfortunately, it is not a priori guaranteed that such a mesh satisfies the constraint on the target minimum dihedral angle $\theta_{\lim}$.
To overcome this difficulty, we apply edge splitting, edge contraction, and point insertion when this stagnation occurs (\cref{pt:coll,pt:split,pt:random} in \cref{alg:allProc}).

For the edge contraction and splitting, we use the standard edge length criterion: we compute the average edge length $l_{\textrm{ave}}$ of the actual mesh, contract the edges shorter than $0.5\,l_{\textrm{ave}}$ (\cref{pt:coll}), and split (halve) the ones longer than $1.5\,l_{\textrm{ave}}$ (\cref{pt:split}).
In \cref{pt:random}, we split a tetrahedron $K$ with $\theta_{\min,K}<\theta_{\lim}$ via a standard 1-to-4 flip by placing the newly added point at the barycenter of $K$~\cite{Edelsbrunner:2006:GTM:1137760}.
In this way, the algorithm constructs via flipping and smoothing a mesh satisfying $Q(\Th)>\theta_{\lim}$.
At the moment, we are not interested in optimizing these operations, we exploit them only to overcome the stagnation of the algorithm.

The MMPDE smoothing can be easily parallelized because the nodal velocities in each smoothing step can be assembled through independent element-wise computation (\cref{eq:vertex:velocity,eq:mesh:eq}), similar to the assembly of a finite element stiffness matrix.
We parallelize the computation of the nodal velocities with OpenMP~\cite{dagum1998openmp}.
Once the velocities are computed, all mesh nodes are moved simultaneously and independently of each other.
On the other hand, the \emph{lazy flip} may propagate to neighbors and neighbors of neighbors, thus, it is complex and difficult to parallelize; in our tests we use a sequential implementation.

\section{Numerical examples}
\label{sec:exe}

We test the proposed mesh improvement algorithm and compare it with the mesh improvement 
algorithm of \Stellar{}~\cite{Klingner:2007:ATM}, the remeshing procedure of 
\CGAL{}~\cite{cgal}, and \MMG{}~\cite{dobrzynski:hal-00681813}.
We compare 
the histograms of the dihedral angles of final meshes, 
the minimal and the maximal dihedral angles $\theta_{\min,\mathcal{T}_h}$ and 
$\theta_{\max,\mathcal{T}_h}$,
the mean dihedral angle $\mu_{\mathcal{T}_h}$, and 
its standard deviation $\sigma_{\mathcal{T}_h}$.

\subsection{Piecewise linear complexes (PLCs)}

\begin{figure*}[p]
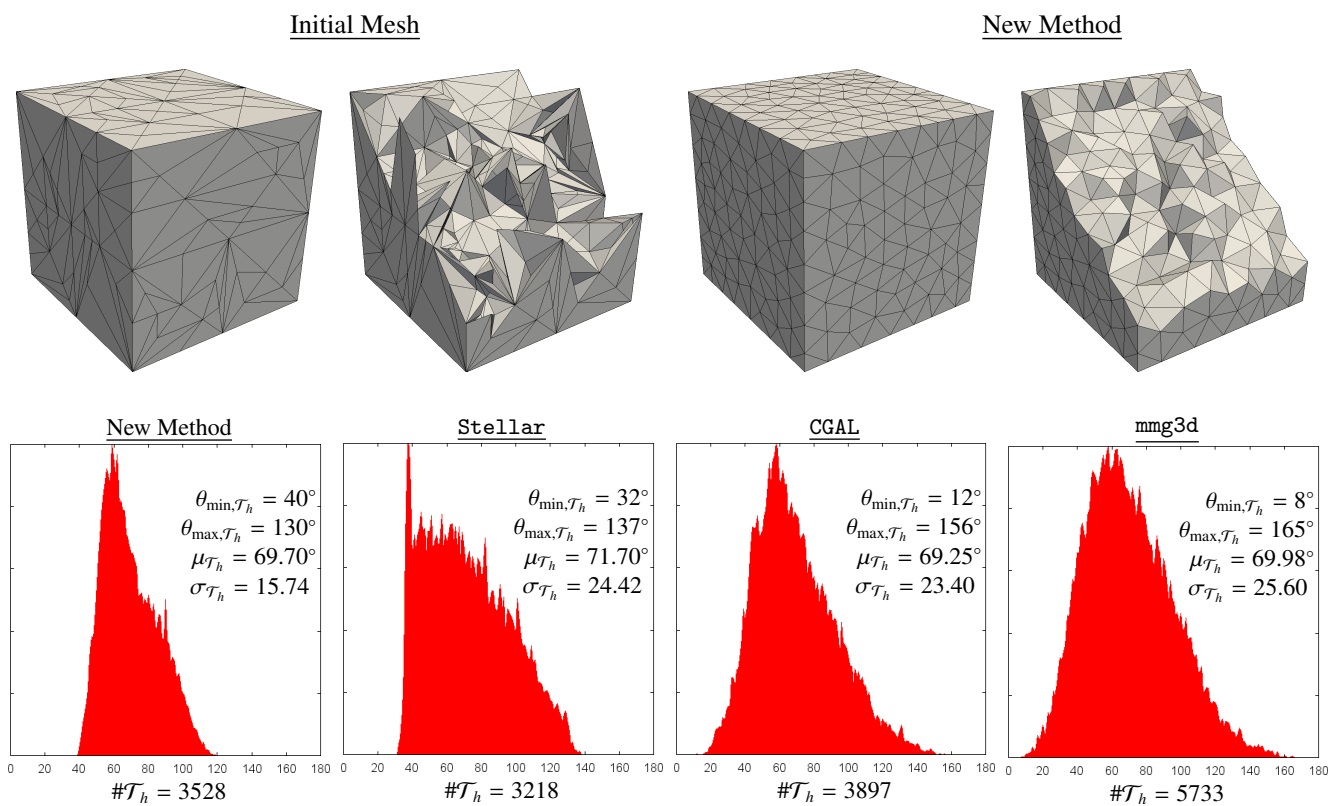

\plotTitle{}
\plotOnlyMesh{Rand1}{ini}{mesh}%
\plotOnlyMesh{Rand1}{ini}{cut}
\plotOnlyMesh{Rand1}{our}{mesh}%
\plotOnlyMesh{Rand1}{our}{cut}\\[-1ex]
\plotResult{New Method}{40}{130}{69.70}{15.74}{806}{3528}{our}{Rand1}%
\plotResult{\Stellar{}}{32}{137}{71.70}{24.42}{1194}{3218}{stellarSameEle}{Rand1}%
\plotResult{\CGAL{}}{12}{156}{69.25}{23.40}{905}{3897}{cgal}{Rand1}%
\plotResult{\MMG{}}{8}{165}{69.98}{25.60}{1086}{5733}{mmg3d}{Rand1}%
\caption{\label{fig:Rand1Hist}
   \textsc{Rand1}.
   The initial mesh with $\#\mathcal{T}_h = \num{5104}$,
   the final (optimized) mesh,
   and statistics of dihedral angles for the final meshes.
}
\end{figure*}

\begin{figure*}[p]
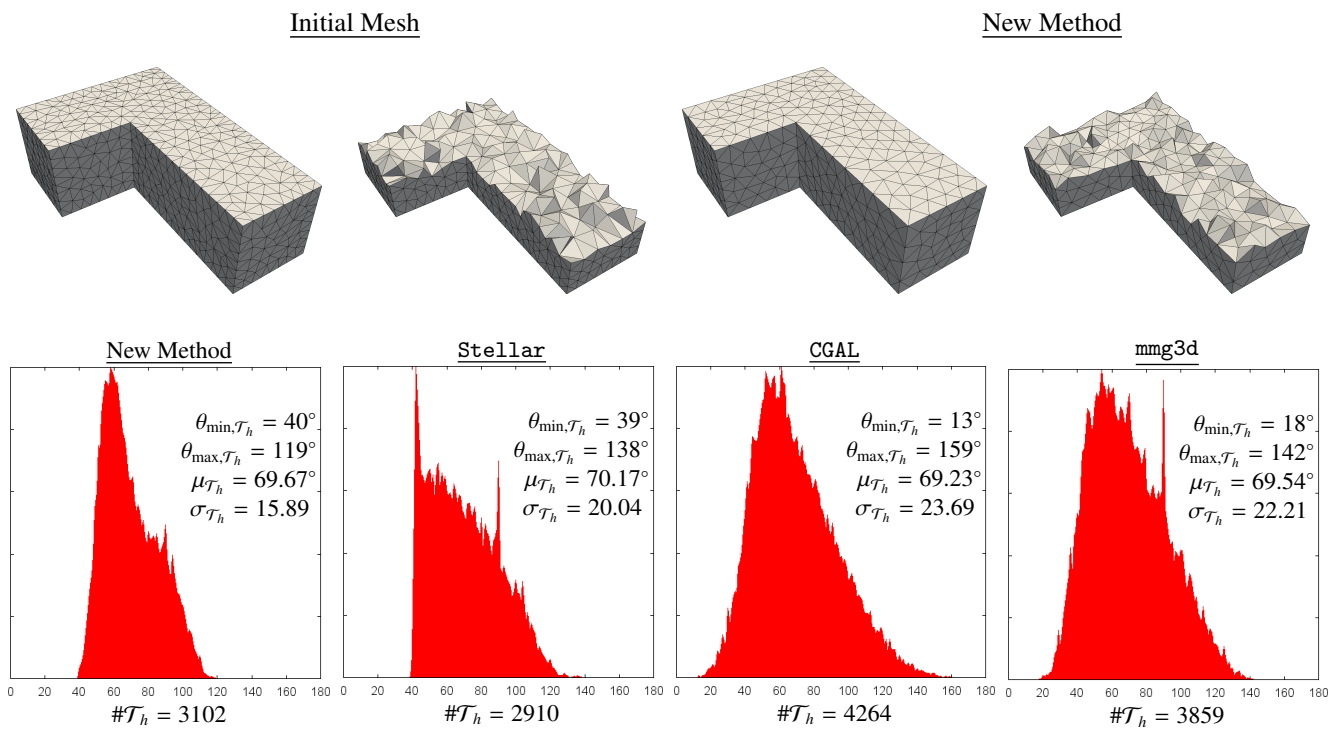

\plotTitle{}
\plotOnlyMesh{LShape}{ini}{mesh}%
\plotOnlyMesh{LShape}{ini}{cut}%
\plotOnlyMesh{LShape}{our}{mesh}%
\plotOnlyMesh{LShape}{our}{cut}\\[-1ex]
\plotResult{New Method}{40}{119}{69.67}{15.89}{796}{3102}{our}{LShape}%
\plotResult{\Stellar{}}{39}{138}{70.17}{20.04}{1258}{2910}{stellar}{LShape}%
\plotResult{\CGAL{}}{13}{159}{69.23}{23.69}{1032}{4264}{cgal}{LShape}%
\plotResult{\MMG{}}{18}{142}{69.54}{22.21}{1114}{3859}{mmg3d}{LShape}%
\caption{\label{fig:LShapeHist}
   \textsc{LShape}.
   The initial mesh with $\#\mathcal{T}_h = \num{4072}$,
   the final (optimized) mesh,
   and statistics of dihedral angles for the final meshes.
}
\end{figure*}

\begin{figure*}[p]
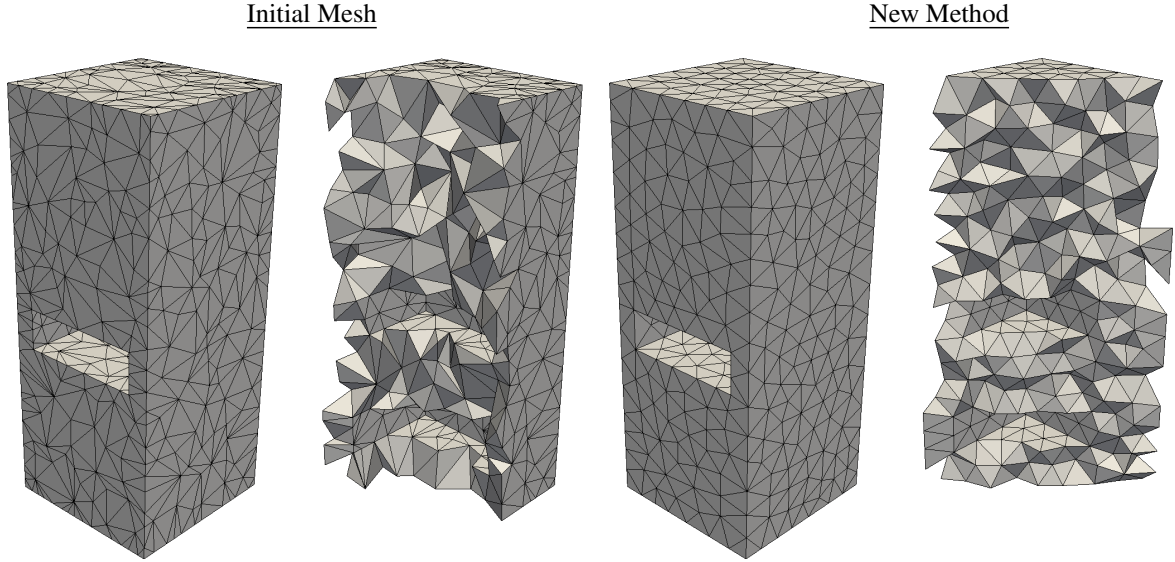
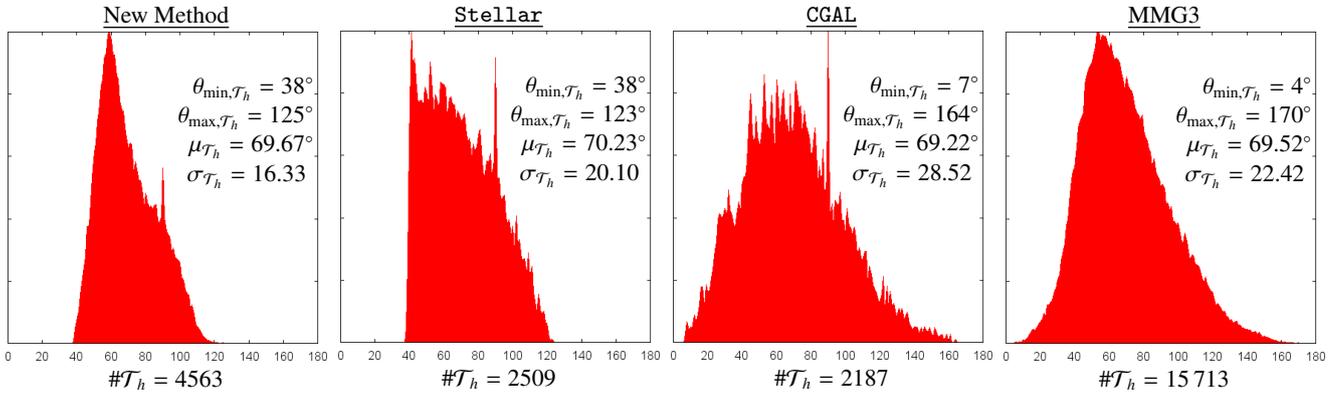
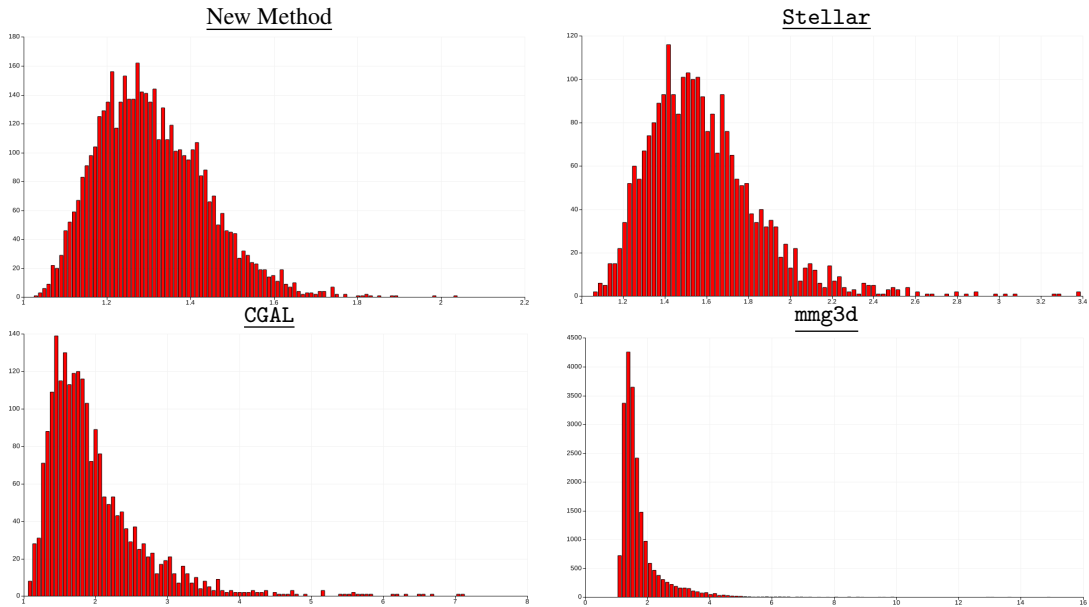
\centering{}
\begin{subfigure}[c]{0.90\textwidth}
\plotTitle{}
\plotOnlyMesh{tetgenExe}{ini}{mesh}%
\plotOnlyMesh{tetgenExe}{ini}{cut}%
\plotOnlyMesh{tetgenExe}{our}{mesh}%
\plotOnlyMesh{tetgenExe}{our}{cut}%
\vspace{-4ex}
\caption{\footnotesize{}
   The initial mesh with $\#\mathcal{T}_h = \num{3545}$
   and the final optimized mesh.\label{fig:tetgen:example:dihedrals}
}
\end{subfigure}
\\[3ex]
\begin{subfigure}[c]{1.0\textwidth}
\plotResult{New Method}{38}{125}{69.67}{16.33}{1120}{4563}{our}{tetgenExe}%
\plotResult{\Stellar{}}{38}{123}{70.23}{20.10}{746}{2509}{stellar}{tetgenExe}%
\plotResult{\CGAL{}}{7}{164}{69.22}{28.52}{660}{2187}{cgal}{tetgenExe}%
\plotResult{MMG3}{4}{170}{69.52}{22.42}{4478}{15713}{mmg3d}{tetgenExe}%
\caption{\footnotesize{}Dihedral angle comparison for the final 
meshes.}\label{fig:tetgenExe:hist}
\end{subfigure}
\\[3ex]
\begin{subfigure}[c]{0.93\textwidth}\centering
\plotAspectRatio{New Method}{our}{tetgenExe}%
\quad
\plotAspectRatio{\Stellar{}}{stellar}{tetgenExe}%
\\
\plotAspectRatio{\CGAL{}}{cgal}{tetgenExe}%
\quad
\plotAspectRatio{\MMG{}}{mmg3d}{tetgenExe}%
\caption{\footnotesize{}Aspect ratio comparison for the final meshes.}
\end{subfigure}
\caption{\label{fig:tetgenExe}
   \textsc{TetgenExample}.
   The initial mesh with $\#\mathcal{T}_h = \num{3545}$,
   the final (optimized) mesh,
   and statistics of the dihedral angles and the aspect ratio.
}
\end{figure*}

To analyze the effectiveness of the proposed mesh improvement scheme 
in case of a piecewise linear complex domain, we consider the following three examples (for more PLC examples, see~\cite{DasKamSi16}):
\begin{itemize}
 \item \textsc{Rand1}  tetrahedral meshes of a cube
 generated by inserting randomly located vertices 
 inside and on the boundary~\cite{Klingner:2007:ATM} (\cref{fig:Rand1Hist}),
 \item \textsc{LShape} is a tetrahedral mesh of an L-shaped PLC
 generated by \texttt{TetGen}~\cite{Si:2015:TDQ:2732672.2629697}
 without optimizing the minimal dihedral angle (switches \texttt{-pa0.019}, 
\cref{fig:LShapeHist}),
 \item \textsc{TetgenExample}
 is a tetrahedral example mesh of a non-convex PLC with a hole
 provided by \texttt{TetGen} (\cref{fig:tetgenExe}).
\end{itemize}

\paragraph{Smoothing and flipping by themselves}
Before testing the full mesh improvement scheme,
we test the effectiveness of the MMPDE smoothing and the \emph{lazy flip}
by themselves and employ smoothing and flipping separately, i.e.,
we improve a tetrahedral mesh exploiting \emph{only} the flipping operation or 
the vertex smoothing.
We compare our results with the ones provided by \Stellar{} for the examples \textsc{LShape} and \textsc{TetGenExample}.

\begin{figure*}[p]
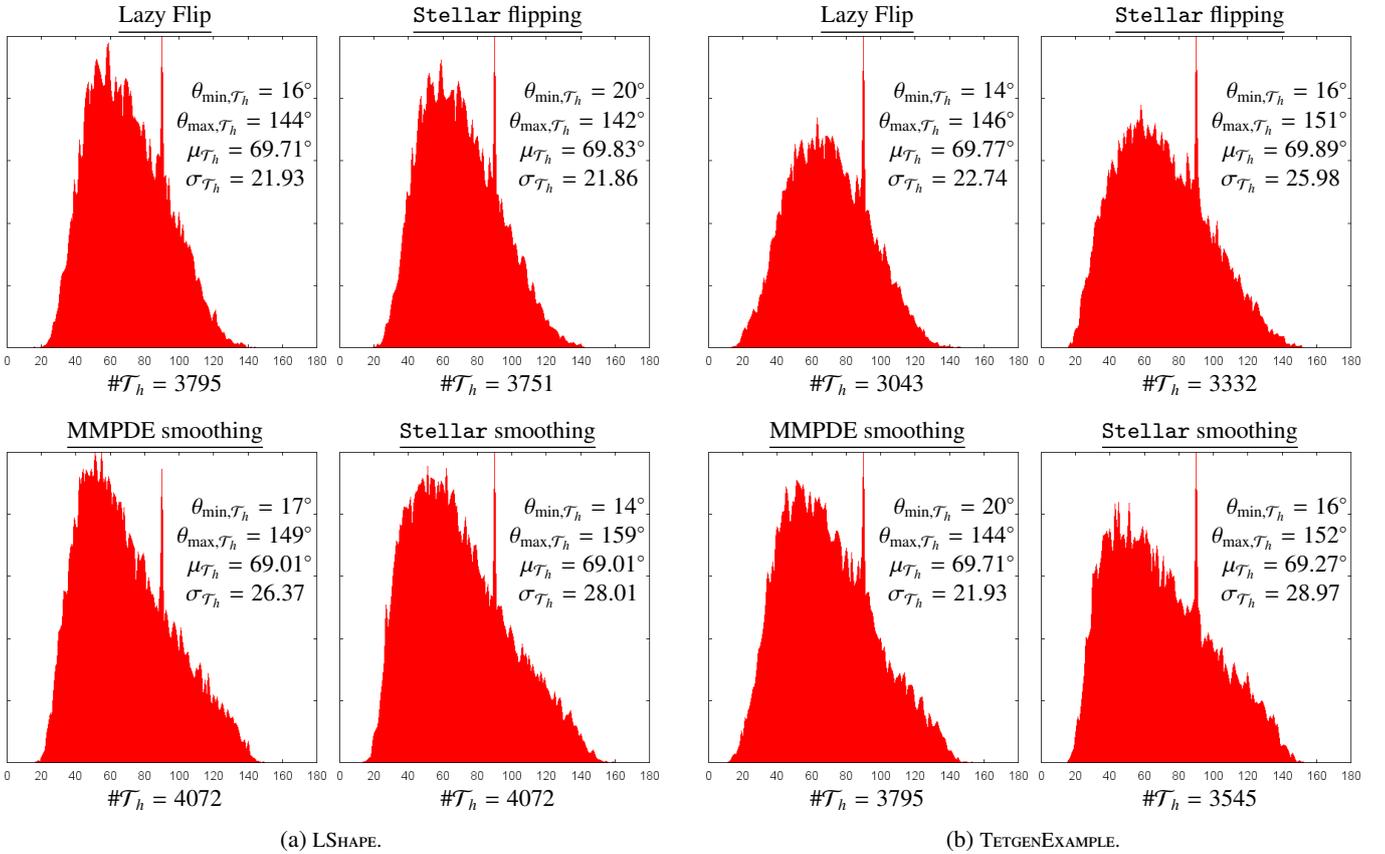

\begin{subfigure}[c]{0.5\textwidth}\centering{}
\plotResultDoubleSize{Lazy Flip} {16}{144}{69.71}{21.93}{1114}{3795}{our}{flip-LShape}%
\plotResultDoubleSize{\Stellar{} flipping} 
{20}{142}{69.83}{21.86}{1}{3751}{stellar}{flip-LShape}
\\[2ex]
\plotResultDoubleSize{MMPDE smoothing} 
{17}{149}{69.01}{26.37}{1114}{4072}{our}{smooth-LShape}%
\plotResultDoubleSize{\Stellar{} smoothing} 
{14}{159}{69.01}{28.01}{1}{4072}{stellar}{smooth-LShape}%
\caption{\footnotesize{}\textsc{LShape}.}\label{fig:smooth:only:lshape}
\end{subfigure}
\begin{subfigure}[c]{0.5\textwidth}\centering{}
\plotResultDoubleSize{Lazy Flip} 
{14}{146}{69.77}{22.74}{881}{3043}{our}{flip-TetgenExample}%
\plotResultDoubleSize{\Stellar{} flipping} 
{16}{151}{69.89}{25.98}{1}{3332}{stellar}{flip-TetgenExample}%
\\[2ex]
\plotResultDoubleSize{MMPDE smoothing} 
{20}{144}{69.71}{21.93}{1114}{3795}{our}{smooth-TetgenExample}%
\plotResultDoubleSize{\Stellar{} smoothing} 
{16}{152}{69.27}{28.97}{1}{3545}{stellar}{smooth-TetgenExample}%
\caption{\footnotesize{}\textsc{TetgenExample}.}\label{fig:smooth:only:tetgen}
\end{subfigure}
\caption{\label{fig:smooth:only}
   Comparison of flipping only (first row) and smoothing only (second row)
   for the initial meshes \textsc{LShape} and \textsc{TetgenExample}.
}
\end{figure*}

The results of the \emph{lazy flip} are comparable to the \Stellar{} flips (\cref{fig:smooth:only}, 
first row).
However, the MMPDE smoothing is better than
its counterpart in \Stellar{} (\cref{fig:smooth:only}, second row):
in both examples it achieves larger $\theta_{\min,\mathcal{T}_h}$,
noticeably smaller $\theta_{\max,\mathcal{T}_h}$,
and a smaller standard deviation of the mean dihedral angle.

\paragraph{Full scheme}
We compare the whole scheme with the mesh improvement algorithm of \Stellar{}~\cite{Klingner:2007:ATM},
the remeshing procedure of \CGAL{}~\cite{cgal}, and \MMG{}~\cite{dobrzynski:hal-00681813}
(\cref{fig:Rand1Hist,fig:LShapeHist,fig:tetgenExe}).

Although all methods provide good results, the new scheme is better: 
$\theta_{\min,\mathcal{T}_h}$ is larger than the value obtained 
by \CGAL{} or \MMG{} and comparable to the value obtained by \Stellar{}.
Moreover, $\theta_{\max,\mathcal{T}_h}$ is smaller than the values obtained by 
\Stellar{}, \CGAL{}, or \MMG{} in all examples but one see \cref{fig:tetgen:example:dihedrals}.

Our method provides mean dihedral angles $\mu_{\mathcal{T}_h}$ around \ang{69.6},
which is close to the optimal value of $\arccos{(1/3)}\approx \ang{70.56}$.
Moreover, standard deviations $\sigma_{\mathcal{T}_h}$ are always smaller than the ones 
of other methods.
Indeed, we get a distribution of dihedral angles close to the mean value.
This quantitative consideration becomes clearer from the shape of the histograms in \cref{fig:Rand1Hist,fig:LShapeHist,fig:tetgenExe}.

For the \textsc{TetgenExample} (\cref{fig:tetgenExe}) we also provide aspect ratio 
histograms (the results for the other examples are very similar and we omit them).
The aspect ratio of an equilateral tetrahedron is equal to $1$ and 
the more a tetrahedron is distorted and stretched the greater its aspect ratio becomes.
Our method and \Stellar{} clearly provide the best aspect ratio distribution. 
For our method, the vast majority of tetrahedra have an aspect ratio smaller than $1.8$.
The \Stellar{} mesh is slightly worse with most of its tetrahedra having aspect ratios below $2.6$.

\subsection{Curved boundary domains}
\label{sec:curved:domains}

In the last part of this section, we experimentally demonstrate some examples with curved domains.
We study two types of examples: 
one academic example for using the RBF surface reconstruction to project the boundary vertices on 
the smooth approximation of the discrete surface and
two more complex examples with fixed boundary vertices.

\begin{figure*}[p]
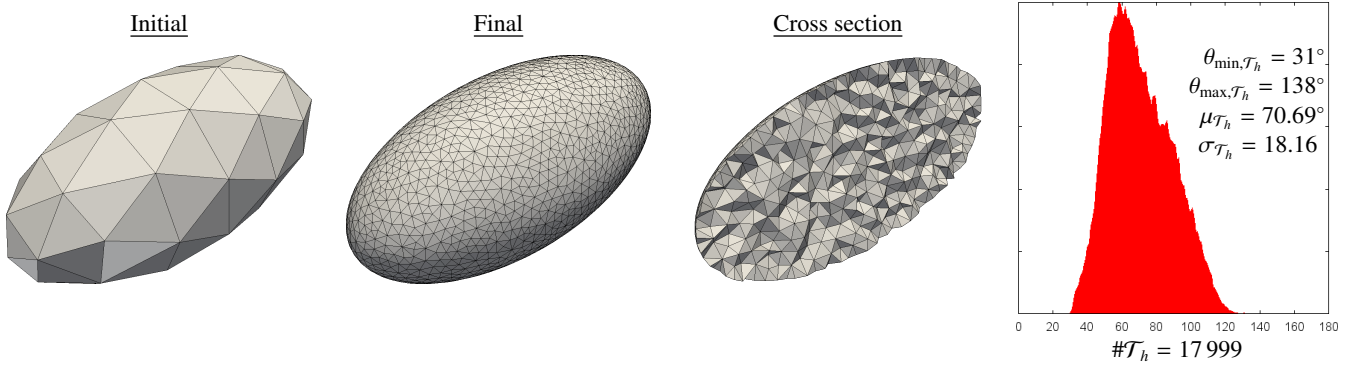

\plotOnlyMeshWithTitle{ellExample}{our}{initial}{Initial}
\plotOnlyMeshWithTitle{ellExample}{our}{final}{Final}
\plotOnlyMeshWithTitle{ellExample}{our}{finalInside}{Cross section}
\plotResult{}{31}{138}{70.69}{18.16}{4122}{17999}{our}{ellExample}
\caption{\footnotesize{} Example of tetrahedral mesh improvement with a curved
(reconstructed) surface.}
\label{fig:ellIni}
\end{figure*}

First, we consider the discrete ellipsoid mesh (\cref{fig:ellIni}).
Though it has a simple geometry, it requires some effort since the boundary is curved and no longer a PLC.\@
The main challenge is to project the boundary vertices back onto the smooth surface if they leave it after a mesh improvement step.
For this reason, we reconstruct the surface via RBFs (see \cref{sub:rfbSurf}) to assist the mesh optimization and project the moved (smoothed) boundary vertices to the reconstructed surface.
\Cref{fig:ellIni} shows that RBF reconstruction smoothes the initially rough surface approximation.
The obtained tetrahedral mesh has high quality: the mean dihedral angle is close to the optimal value ($\approx\ang{70.69}$) and the standard deviation of the dihedral angles is small ($\approx\ang{18.16}$).

However, it has to be pointed out that complicated boundaries cannot be handled as easily as an ellipsoid and require more sophisticated methods.

In our next examples, we restrict ourselves to the case of fixed boundary vertices since \Stellar{} does not handle curved surfaces described via an implicit function, start with a good isotropic triangular mesh as input, and keep the boundary vertices fixed for each of the algorithms.

\begin{figure*}[p]
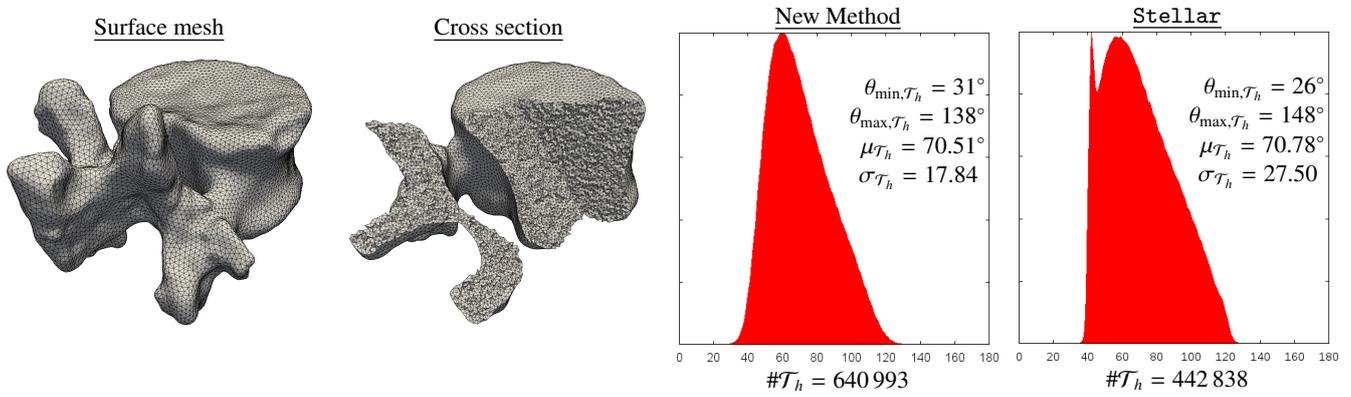

\plotOnlyMeshWithTitle{spineExe}{our}{spineAll}{Surface mesh}
\plotOnlyMeshWithTitle{spineExe}{our}{spineSec}{Cross section}
\plotResult{New Method}{31}{138}{70.51}{17.84}{4122}{640993}{our}{spineExe}
\plotResult{\Stellar{}}{26}{148}{70.78}{27.50}{4122}{442838}{stellar}{spineExe}
\caption{\footnotesize{} Spine example: the initial mesh with $\#\mathcal{T}_h = 
\num{688420}$ and the final optimized mesh.}
\label{fig:spineRes}
\end{figure*}

\begin{figure*}[p]
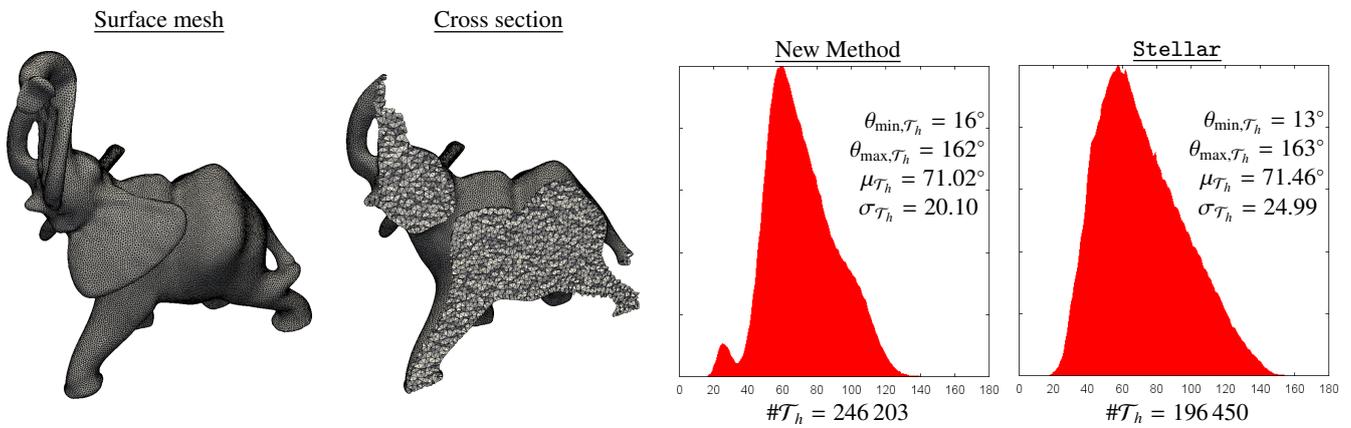

\plotOnlyMeshWithTitle{elephantExe}{our}{meshAll}{Surface mesh}
\plotOnlyMeshWithTitle{elephantExe}{our}{meshInside}{Cross section}
\plotResult{New Method}{16}{162}{71.02}{20.10}{4122}{246203}{our}{elephantExe}
\plotResult{\Stellar{}}{13}{163}{71.46}{24.99}{4122}{196450}{stellar}{elephantExe}
\caption{\footnotesize{} Elephant example: the initial mesh with $\#\mathcal{T}_h = 
\num{260401}$ and the final optimized mesh.}
\label{fig:elephantRes}
\end{figure*}

\paragraph{Fixed boundary}
The next two examples are meshes of a spinal bone and of an elephant (\cref{fig:spineRes,fig:elephantRes}).
The initial surface meshes in both examples are constructed by means of the higher dimensional embedding approach for surface mesh reconstruction~\cite{DasFarSi16} and their minimal face angles are approximately \ang{33}.
The initial volume meshes are constructed by \texttt{TetGen} using the \texttt{-Y} flag to preserve the fixed boundary.

\Cref{fig:spineRes,fig:elephantRes} present the histograms of the dihedral angles of the resulting optimized meshes.
In comparison to the PLC examples, where the geometry is simpler and the boundary vertices are allowed to move, the smallest dihedral angles for the spinal bone and the elephant examples are worse (smaller) than for the PLC examples.
In comparison to \Stellar{}, our algorithm achieves better values for $\theta_{\min,\mathcal{T}_h}$ and $\theta_{\max,\mathcal{T}_h}$, as well as a smaller mean deviation from the mean value.

These examples, too, show the ``aggressive'' nature of the \Stellar{} mesh improvement algorithm, which aggressively removes vertices during the mesh improvement.
In contrast, our mesh improvement scheme is able to produce a high-quality mesh while keeping the number of vertices close to the original input.

\section{Conclusions}\label{sec:conc}

Mesh improvement is a challenging problem and we tackled it by combining several recently developed techniques, namely, moving mesh smoothing, lazy flipping, and RBF surface reconstruction.
In comparison to the mesh improvement algorithm \Stellar{} and the re-meshing procedures provided by \CGAL{} and \MMG{}, 
we obtain better results in terms of the distributions of dihedral angles for all 
 considered examples.
However, there are several directions in which this work could be extended.

First, for smooth and relatively simple boundaries, our approach works excellently but complicated curved boundaries pose a challenging problem.
One possible solution could be the direct incorporation of the boundary description into the MMPDE smoothing scheme (parametrization) so that the boundary vertices will always stay on the surface.
This will avoid the sometimes troublesome projection of vertices and velocities back onto the surface after a smoothing step.

Second, we need to find a more sophisticated method for edge contraction and splitting in order to improve the performance of both the MMPDE smoothing and the lazy flip.

Third, the MMPDE smoothing is based on the moving mesh method~\cite{Huang2015322} which allows the definition of a metric field.
Hence, the moving mesh smoothing can be extended to the adaptive and anisotropic setting.


\section*{Acknowledgments}
The work of Franco Dassi was partially supported by the ``Leibniz-DAAD Research Fellowship 2014''.
The authors are thankful to Jeanne Pellerin for her support in computing the examples with \MMG{}~\cite{dobrzynski:hal-00681813}
and the anonymous referee for the valuable comments
which helped to improve the quality of this paper.



\bibliographystyle{elsarticle-num}
\bibliography{DasKamFarSi17}

\end{document}